\documentclass{elsarticle}
\usepackage{amssymb, amsmath, amsthm}

\newtheorem{theorem}{Theorem}[section]
\newtheorem{lemma}[theorem]{Lemma}
\newtheorem{proposition}[theorem]{Proposition}
\newtheorem{corollary}[theorem]{Corollary}

\newtheorem{definition}[theorem]{Definition}

\newtheorem{remark}[theorem]{Remark}

\numberwithin{equation}{section}

\newcommand\NN{\mathbb{{N}}}

\newcommand\RR{\mathbb{{R}}}
\newcommand\CC{\mathbb{{C}}}
\newcommand\ZZ{\mathbb{{Z}}}

\newcommand\ra{{\rm a}}

\newcommand\rc{{\rm c}}
\newcommand\rd{{\rm d}}

\newcommand\ba{{\bf a}}
\newcommand\bb{{\bf b}}

\newcommand\bd{{\bf d}}

\newcommand\bbf{{\bf f}}

\newcommand\br{{\bf r}}
\newcommand\bt{{\bf t}}
\newcommand\bv{{\bf v}}

\newcommand\bx{{\boldsymbol x}}
\newcommand\by{{\bf y}}
\newcommand\bz{{\boldsymbol z}}

\newcommand\bT{{\bf T}}

\newcommand\balpha{{\boldsymbol\alpha}}
\newcommand\bbeta{{\boldsymbol\beta}}
\newcommand\bdelta{{\boldsymbol \delta}}
\newcommand\bfeta{{\bf \eta}}
\newcommand\bepsilon{{\boldsymbol\epsilon}}

\newcommand\bgamma{{\boldsymbol\gamma}}
\newcommand\bomega{{\boldsymbol\omega}}
\newcommand\blambda{{\boldsymbol\lambda}}

\newcommand\btau{{\boldsymbol\tau}}
\newcommand\bvarepsilon{{\boldsymbol\varepsilon}}
\newcommand\bxi{{\boldsymbol\xi}}

\newcommand\bnull{{\boldsymbol0}}
\newcommand\bone{{\boldsymbol1}}

\newcommand{\mfrac}[2]%
{\raisebox{0.5pt}{\footnotesize$\dfrac{#1}{#2}$}}
\newcommand{\mbinom}[2]%
{\raisebox{0.5pt}{\footnotesize$\dbinom{#1}{#2}$}}
\def\smmat\{#1&#2\cr#3&#4\}%
   {\!\!\pmatrix{\!#1\!&\!#2\!\cr\!#3\!&\!#4\!}\!\!}
\newcommand\scrm{{\raise0.5pt\hbox{-}}}
\def\eop{{ \vrule height7pt width7pt depth0pt}\par\bigskip} 
\newcommand\ie{{\it\thinspace i.e.}}
\newenvironment{pf}%
{\par\noindent\textbf{Proof:}~}

\begin{document}

\title{Reproduction of Exponential Polynomials by Multivariate Non-stationary
Subdivision Schemes with a General Dilation Matrix}
\author[ch]{Maria Charina\corref{cor1}}
\ead{maria.charina@uni-dortmund.de}
\author[cc]{Costanza Conti}
\ead{costanza.conti@unifi.it}
\author[lr]{Lucia Romani}
\ead{lucia.romani@unimib.it}

\cortext[cor1]{Corresponding author}
\address[ch]{Fakult\"at f\"ur Mathematik, TU Dortmund, D--44221 Dortmund, Germany}
\address[cc]{Dipartimento di Energetica, Universit\`a di Firenze, Via C. Lombroso 6/17, I--50134 Firenze}
\address[lr]{Dipartimento di Matematica e Applicazioni, Universit\`{a} di Milano-Bicocca, Via R. Cozzi 53, 20125 Milano, Italy}

\begin{abstract}
 We study scalar multivariate non-stationary subdivision schemes with a general dilation matrix.
 We characterize the capability of such schemes to reproduce exponential polynomials in terms of simple
 algebraic conditions on their symbols. These algebraic
 conditions provide a useful theoretical tool for checking the
 reproduction properties of existing schemes and for
 constructing new schemes with desired reproduction capabilities and other enhanced properties.
We illustrate our results with several examples.
\end{abstract}


\maketitle

\section{Introduction}

In this paper we study multivariate non-stationary subdivision schemes which are iterative algorithms based on the repeated application of
subdivision operators \begin{equation} \label{def:sub_operator}
 S_{\ba^{[k]}}: \ell(\ZZ^s) \rightarrow \ell(\ZZ^s), \quad \bbf^{[k+1]}= S_{\ba^{[k]}} \bbf^{[k]}, \quad k \ge 0.
\end{equation}
The linear subdivision operators $ S_{\ba^{[k]}}$
$$
 \left(S_{\ba^{[k]}} \bbf^{[k]}\right)_\balpha=\sum_{\bbeta \in \ZZ^s} \ra^{[k]}_{\balpha-M\bbeta}
f^{[k]}_\bbeta,\qquad \balpha \in\ZZ^s,
  \qquad k \ge 0.
$$
depend on a dilation matrix $M \in \ZZ^{s \times s}$ and on the finite sequences of real
numbers $\ba^{[k]}=\{ \ra_\balpha^{[k]}, \ \balpha\in \ZZ^s\}$.
If convergent, subdivision schemes are used e.g. for designing  curves, surfaces, or multivariate functions.

\smallskip \noindent
The \emph{generation} properties of subdivision schemes  are well-understood and are
characterized in terms of so-called \emph{zero conditions}, see e.g. \cite{CDM,DLL03,JePlo},
on the mask symbols
$$
 a^{[k]}(\bz)=\sum_{\balpha \in \ZZ^s} \ra^{[k]}_\balpha \bz^\balpha, \quad k \ge
 0 \quad \bz \in \left( \CC \setminus\{0\} \right)^s.
$$
These conditions on the symbols determine if the subdivision limit belongs to the
function space 
$$
 EP_{\Gamma, \Lambda}=\hbox{span}\{\bx^\bgamma e^{\blambda \cdot \bx} : \ \bgamma \in \Gamma, \ \blambda \in \Lambda
 \}, \quad \Gamma \subset \NN_0^s, \quad \Lambda \subset \CC^s,
$$
if the starting sequence $\bbf^{[0]}$ in \eqref{def:sub_operator} is sampled from  a function in this
space. If for all functions in  $EP_{\Gamma, \Lambda}$, the subdivision scheme
generates exactly the same function from which the corresponding starting sequence $\bbf^{[0]}$ is
sampled, then we say that the scheme \emph{reproduces} $EP_{\Gamma, \Lambda}$.

\smallskip \noindent
The main goal of our study is to characterize the reproduction property of subdivision
in terms of zero conditions and some additional conditions on  subdivision symbols.
These additional  algebraic conditions
provide a useful theoretical tool that simplifies the analysis of reproduction properties
of existing schemes and the construction of new schemes with desired reproduction properties.
We would like to emphasize that the properties of non-stationary subdivision schemes could be 
characterized in the Fourier domain in terms of the properties of the associated basic refinable functions, see e.g 
\cite{BR, DR, Han12, JiaLei, Ron_exp_box_splines, Ron92}. Nevertheless,  we
think that it is more advantageous, for construction of new non-stationary schemes, to provide 
such characterizations in terms of  simple algebraic conditions on  the coefficients
$\ba^{[k]}$, $k \ge 0$, of the corresponding refinement equations. 

\smallskip \noindent
The reproduction and generation of $EP_{\Gamma, \Lambda}$ or its subspaces are crucial for modeling objects of different
shapes that e.g. are described by polynomial, trigonometric or hyperbolic functions.
Thus, these reproduction and generation properties of subdivision are
important in CAGD, motion planning, iso-geometric analysis and for studying approximation properties
of subdivision schemes \cite{ALevin, Yoon}. In spite of their importance, reproduction properties
of non-stationary subdivision schemes have not been yet studied as rigorously as it is done in the stationary case when
the mask $\ba$ does not depend on the level $k$ of subdivision recursion. In the stationary case, characterizations of
reproduction of polynomial spaces (corresponding to $\Lambda=\{(0,\ldots,0)\}$) are given in
\cite{CC2012, ContiHormann11, Han03, Han13} where, in addition to zero conditions on the symbol $a(\bz)$, extra conditions on
$a(\bz)$ and on its derivatives at $(1,\ldots,1)$ are established.

\smallskip \noindent Characterizations
of reproduction properties of binary univariate non-stationary schemes are given
in \cite{CR2011} and we extend them here to the case of multivariate non-stationary subdivision
 with a general integer dilation matrix $M$. We study the so-called non-singular schemes, i.e
 the ones that generate a zero function if and only if $\bbf^{[0]}$ is a zero sequence.
The main result of our paper states
that a non-stationary non-singular subdivision scheme defined by the masks $\{\ba^{[k]}, \ k \ge 0\}$
reproduces $EP_{\Gamma, \Lambda}$
if and only if there exists $\btau \in \RR^s$ such that the symbols $a^{[k]}(\bz)$ and their
derivatives $D^\bgamma a^{[k]}(\bz)$ satisfy the zero conditions \eqref{algebraic_conditions_multiGENE} and 
$$
 \bv^\bgamma D^\bgamma a^{[k]}(\bv)= |\hbox{det}(M)| \cdot \bv^{M\btau-\btau}  \prod_{\ell=1}^s \prod_{j=0}^{\gamma_\ell-1}
  ((M\btau-\btau)_\ell-j),
$$
for all $\bv=e^{-\blambda \cdot M^{-(k+1)}}$, $\bgamma \in \Gamma$, $\blambda \in \Lambda$ and $k\ge 0$, see
Theorem \ref{th:characterization_poly_reproduction_multi} for details.

\smallskip \noindent For a convergent non-stationary subdivision scheme, these
conditions are sufficient for reproduction of $EP_{\Gamma, \Lambda}$,
see subsection \ref{subsec:reproduction_multi}. The
parameter $\btau \in \RR^s$ above is the so-called \emph{shift} parameter that determines the
\emph{parametrization}
\begin{equation} \label{def:parametrization}
 \bt^{[k]}_{\balpha}=M^{-k}(\balpha+\btau),\quad  \balpha \in \ZZ^s, \quad \btau \in \RR^s,
 \quad k\ge 0,
\end{equation}
associated to a given non-stationary subdivision scheme. This
parametrization specifies to which grid points $\bt^{[k]}_{\balpha}$ the newly computed values
$f_\balpha^{[k]}$
are attached at the $k-$th level of subdivision recursion. The parametrization also influences the starting sequences $\bbf^{[0]}$.
The choice of $\btau$ does not affect either the generation property or convergence of a
subdivision scheme, but it affects its reproduction properties. In particular, a  correctly chosen shift parameter allows to enrich
the variety of shapes the scheme can reproduce.

%

\smallskip \noindent The paper is organized as follows: in Section \ref{sec:notation}, we introduce notation and
recall some known facts about non-stationary subdivision. In Section \ref{sec:new}, for a non-stationary subdivision scheme,
we define the notions of generation and reproduction of spaces of exponential polynomials and 
emphasize the difference between these notions. Since the univariate case is
certainly simpler to follow, in Section  \ref{sec:algebra}, we first derive the algebraic
conditions that characterize the reproduction of exponential polynomials in the
univariate $m$-ary case, $M=m$, $m \ge 2$. Then, in Subsection  \ref{subsec:multi_generation_reproduction}, we extend this characterization to
the multivariate setting. Finally, in Section  \ref{sec:examples}, we apply the derived algebraic
conditions and study subdivision schemes of any arity associated to exponential
B-splines and exponential box splines. We also show the effect of the correct
choice of $\btau$ and also the effect of the renormalization of subdivision
symbols on the reproduction properties of subdivision. We conclude Section 5
with several univariate and bivariate examples. All examples illustrate that our
algebraic conditions on subdivision symbols make the construction of new
schemes with desired reproduction capabilities and with other enhanced
 properties more efficient.

\section{Notation and Subdivision Background} \label{sec:notation}
We start this section by setting the notation and continue by  recalling some known facts about
non-stationary subdivision schemes.

\begin{itemize}
\item[-] $\NN_0$ is the set of natural numbers that includes  zero;
\item[-] Multi-indices are denoted by Greek boldface letters
$\balpha=(\alpha_1,\ldots,\alpha_s)^T \in \NN_0^s$;
\item[-] Vectors are denoted by boldface letters  $\bx \in \RR^s$, or $\bz\in\CC^s$;
\item[-] In the multi-index notation we have $\bx^{\alpha}=x_1^{\alpha_1} \cdot  \ldots  \cdot  x_s^{\alpha_s}$,
$\balpha!=\alpha_1!\ldots\alpha_s!$ and $|\balpha|=\alpha_1+ \ldots + \alpha_s$ for $\balpha\in \NN_0^s$ and $\bx\in \RR^s$;
\item[-]The product of a real number $y\in \RR$ and a vector $\bx$ is denoted by $y \,
\bx=(y x_1, \ldots ,y x_s)^T$;
\item[-] The scalar product of two $s$-dimensional (column) vectors $\bx$ and $\bz$ is
denoted by $\bx\cdot  \bz=x_1z_1  + \ldots + x_sz_s$;
\item[-] The  scalar product $\bx^T M$ of an $s$-dimensional (column) vector $\bx$ and a matrix
 $M$ is denoted by $\bx \cdot M$;
\item[-] For vectors $\bx$ and $\by$, we have $e^\bx=(e^{x_1},\cdots,e^{x_s})^T$ and $(e^\bx)^\by=e^{\bx\cdot \by}$;
\item[-] Fourier transform $\hat{f}$ of a compactly supported essentially bounded function $f \in L_\infty(\RR^s)$
  is defined by  $\displaystyle \hat{f}(\bomega)=\int_{\RR^s} f(\bx)e^{-i \bx \cdot \bomega}d\bx, \ \ \bomega \in \RR^s$;
\item[-] A sequence of real numbers indexed by $\ZZ^s$ is denoted by boldface letters
    $\ba=\{\ra_\balpha, \  \balpha\in \ZZ^s\}$;
\item[-] A sequence of $s$-dimension vectors indexed by $\ZZ^s$ is denoted by capital
        boldface letters $\bT=\{\bt_\balpha, \ \balpha\in \ZZ^s\}$;
\item[-] The space of bounded sequences indexed by $\ZZ^s$ is denoted by $l_\infty(\ZZ^s)$;
\item[-] By $D^\bgamma$, $\bgamma \in \NN_0^s$, we denote a directional derivative.
\end{itemize}

\subsection{Non-stationary subdivision schemes}

A \emph{non-stationary} multivariate \emph{subdivision scheme} is an iterative algorithm
with refinement rules
\begin{equation}\label{def:Soperator0}
  \left(S_{ \ba^{[k]}} \bbf^{[k]}\right)_\balpha=f^{[k+1]}_\balpha = \sum_{\bbeta \in \ZZ^s} \ra^{[k]}_{\balpha-M\bbeta}
f^{[k]}_\bbeta,\qquad \balpha \in\ZZ^s,
  \qquad k \ge 0 \,,
\end{equation}
and it generates the refined data sequence $\bbf^{[k+1]} = \{f^{[k+1]}_\balpha , \  \balpha \in \ZZ^s\}$  from $\bbf^{[k]} = \{f^{[k]}_\balpha , \  \balpha \in \ZZ^s\}$. The dilation matrix $M \in \ZZ^{s \times s}$
is assumed to have all its eigenvalues greater than $1$ in absolute value.
Subdivision schemes are
based on the application of the subdivision operators $S_{\ba^{[k]}}$ constructed from
the so-called \emph{subdivision mask} $\ba^{[k]}=\{\ra^{[k]}_\balpha   , \  \balpha \in \ZZ^s \}$ at level $k$. 
Each  subdivision mask is assumed to be of finite support,
\ie\  $supp(\ba^{[k]})\subset[0,N_k]^s$ for some positive integer $N_k$.
The (non-stationary) subdivision scheme is denoted  by
$S_{\{\ba^{[k]}, \  k\ge 0\}}$ and is given by

$$
\left\{
  \begin{array}{ll}
    \hbox{Input}:\quad \bbf^{[0]},\  \{\ba^{[k]}, \  k\ge 0\}\\
\\
    \hbox{For}\ k=0,1,2,\ldots,\\
\quad \quad \bbf^{[k+1]}=S_{\ba^{[k]}}\bbf^{[k]}
  \end{array}
\right.
$$
for some initial data $\bbf^{[0]}=\{f^{[0]}_\balpha , \  \balpha \in \ZZ^s \}$.

\smallskip \noindent
The \emph{symbols} of a non-stationary subdivision scheme are given by  Laurent polynomials
\begin{equation}\label{eq:laurpol}
  a^{[k]}(\bz) = \sum_{\balpha\in\ZZ^s} \ra^{[k]}_\balpha\;\bz^\balpha, \quad
  \bz \in\left( \CC\setminus\{0\}\right)^s,\quad k\ge 0\,.
\end{equation}
Denote  $m=|\hbox{det}(M)|$ and  by
$E$  the set of representatives of $\ZZ^s / M \ZZ^s$. Clearly, $E$ contains
$\bnull=(0,\ldots,0)^T$. Define the  set
\begin{equation}\label{def:Xi}
  \Xi =\{ e^{2\pi  i  M^{-T}\bxi} \ : \   \bxi \ \hbox{is a coset representative of} \ \ZZ^s / M^T \ZZ^s\}\,,
\end{equation}
which contains $\bone=(1, 1, \dots, 1)^T$. The $m$ \emph{sub-masks} of the masks $ \ba^{[k]}$ are
\begin{equation} \label{eq:submasks}
   \{\ra^{[k]}_{\bvarepsilon+M\balpha},\ \balpha\in\ZZ^s\},\quad \bvarepsilon\in E,\quad k\ge 0\,,
  \end{equation}
and  their symbols (\emph{sub-symbols }of the masks) are
  \begin{equation} \label{eq:subsymbols}
  \quad a^{[k]}_\bvarepsilon(\bz) = \sum_{\balpha\in\ZZ^s}
  \ra^{[k]}_{\bvarepsilon+M\balpha}\;\bz^{\bvarepsilon+M\balpha}
  \;,\quad \bvarepsilon\in E,\quad k\ge 0\,,
\end{equation}
respectively.
The symbols satisfy
\begin{equation} \label{eq:decomp}
  a^{[k]}(\bz) =
  \sum_{\bvarepsilon\in E} \; a^{[k]}_{\bvarepsilon}(\bz),\quad k\ge 0\,.
\end{equation}
For a given parametrization
$\bT^{[k]}=\{\bt_\balpha^{[k]},\ \balpha\in \ZZ^s\}$ 
in \eqref{def:parametrization}, a notion of convergence for  $S_{\{\ba^{[k]}, \ k\ge 0\}}$ is
established using
the sequence $\{F^{[k]}, \  k \ge 0 \}$ of continuous functions $F^{[k]}$
that interpolate the data $\bbf^{[k]}$ at the parameter values $\bt_\balpha^{[k]}$, $\balpha\in \ZZ^s$,
namely
\begin{equation} \label{eq:interpolating_functions}
  F^{[k]}( \bt_\balpha^{[k]}) = f_\balpha^{[k]}, \qquad \balpha \in \ZZ^s,\quad  k\geq0.
\end{equation}
The scheme $S_{\{\ba^{[k]}, \  k\ge 0\}}$ applied to initial data $\bbf^{[0]}$ is called \emph{convergent}, if there
exists a continuous \emph{limit function} $g_{\bbf^{[0]}}$ (which is nonzero for at least one initial nonzero
sequence) such that the sequence $\{F^{[k]}, \ k \ge 0 \}$ converges
uniformly to $g_{\bbf^{[0]}}$, i.e.
\begin{equation}\label{eq:subdivisionlimit}
  g_{\bbf^{[0]}} =\lim_{k \to\infty} S_{\ba^{[k]}} S_{\ba^{[k-1]}} \cdots S_{\ba^{[0]}} \bbf^{[0]}= \lim_{k \to\infty} F^{[k]}.
\end{equation}
The scheme
$S_{\{\ba^{[k]}, \ k \ge 0\}}$  is called \emph{weakly convergent},  if the sequence
$\{F^{[k]}, \ k \ge 0 \}$ in \eqref{eq:interpolating_functions} converges
pointwise to $g_{\bbf^{[0]}} \in L_\infty(\RR^s)$.

\smallskip \noindent
Due to linearity of each subdivision operator $S_{\ba^{[k]}}$, the limit functions
$g_{\bbf^{[0]}}$
exist if and only if the subdivision scheme applied
to the initial data $\bdelta=\{\delta_{\balpha, \bnull} , \  \balpha\in\ZZ^s\}$
converges to the so-called \emph{basic
limit function} $\phi=g_{\bdelta}$. In this case,
$$
 g_{\bbf^{[0]}}=\sum_{\balpha\in \ZZ^s}f^{[0]}_\balpha \phi(\cdot-\balpha).
$$

\smallskip \noindent Differently from the stationary case where all masks
$\ba^{[k]}$ are the same, in the non-stationary setting, one could start the
subdivision process with a mask at level $\ell \geq 0$ and get a family of
subdivision schemes based on the masks $\{\ba^{[\ell+k]}, \  k\ge 0\},\ \ell \geq 0$.
The corresponding subdivision limits are denoted by
\[
  g^{[\ell]}_{\bbf} := \lim_{k\to\infty} S_{\ba^{[\ell+k]}}S_{\ba^{[\ell+k-1]}} \cdots
S_{\ba^{[\ell]}} \bbf^{[0]},\quad \ell \ge 0\,.
\]
An interesting fact about the compactly supported basic limit functions
$\phi^{[\ell]}:=g^{[\ell]}_{\bdelta},\ \ell \geq 0$, ($\phi^{[0]}=\phi$) is that they are mutually refinable, i.e.,
they satisfy the functional equations
\begin{equation}\label{eq:ns-refinability}
    \phi^{[\ell]}=\sum_{\balpha \in \ZZ^s} \ra^{[\ell]}_\balpha \phi^{[\ell+1]}(M\cdot-\balpha),\quad \ell \geq 0,\quad
\hbox{with}\ \ \ba^{[\ell]}  \ \ \hbox{the $\ell$-th level mask}.
\end{equation}

In this paper, we  consider subdivision schemes that are \emph{non-singular}, i.e. they are convergent
 (or weakly convergent) and, such that $g_{\bbf^{[0]}}=0$ if and
only if $\bbf^{[0]}=0$. The following proposition connects the non-singularity of a non-stationary scheme with the linear
independence of the
translates of $\phi^{[\ell]}$, $\ell \ge 0$.

\begin{proposition} \label{prop:nonsingular_vs_independence}
If a subdivision scheme $S_{\{\ba^{[k]},\ k\ge 0\}}$  is non-singular, then the integer translates of
$\phi^{[\ell]}$ are linearly independent for each $\ell \ge 0$.
\end{proposition}
\pf
\noindent
If the convergent subdivision scheme $S_{\{\ba^{[k]}, \ k\ge 0\}}$ is non-singular, then
for any starting sequence $\bd=\{\rd_\balpha, \ \balpha \in \ZZ^s\}$ we have
$$
 \lim_{k \rightarrow \infty} S_{\ba^{[k+\ell]}} \dots S_{\ba^{[\ell]}}  \bd=\sum_{\balpha \in \ZZ^s} \rd_\balpha \phi^{[\ell]}(\cdot-\balpha)=0
$$
if and only if  $\bd$ is the zero sequence.
\eop

We formulate the converse of Proposition \ref{prop:nonsingular_vs_independence} separately
as it requires  additional assumptions on the functions $\phi^{[\ell]}$.

\begin{proposition} \label{prop:independence_vs_nonsingular}
Let $\phi^{[\ell]} \in L_\infty(\RR^s)$, $\ell \ge 0$, be compactly supported
solutions of the refinement equations \eqref{eq:ns-refinability} and such that
their integer shifts are linearly independent for each $\ell$.  If there exists
$L \ge 0$ such that $\phi^{[\ell]}$ additionally satisfy
$$
 C \| \bbf \|_\infty \le \|\sum_{\balpha \in \ZZ^s} f_\balpha \phi^{[\ell]}(\cdot-\balpha)\|_\infty,
 \quad \bbf \in l_\infty(\ZZ^s), \quad \ell \ge L,
$$
then the scheme $S_{\{\ba^{[k]}, \ k\ge 0\}}$ is non-singular.
\end{proposition}
\pf
 By \cite[Theorem 13]{DL95}, the assumptions on $\phi^{[\ell]}$ guarantee that the associated
 subdivision scheme $S_{\{\ba^{[k]},\ k\ge 0\}}$ is convergent. The linear independence of
 $\phi^{[\ell]}$ yields the claim.
\eop

\section{Exponential polynomials and non-stationary subdivision schemes} \label{sec:new}

In this section, for a non-stationary subdivision scheme, we define the notions of
\emph{exponential polynomial generation} and \emph{exponential polynomial reproduction}.
We start by
defining the space of \emph{exponential polynomials} on $\RR^s$.

\begin{definition} For $\Gamma \subset \NN_0^s$ and $\Lambda \subset \CC^s $  we define
$$
 EP_{\Gamma,\Lambda}=\hbox{span}\{\bx^{\bgamma} e^{\blambda \cdot \bx}, \  \ \bgamma \in \Gamma, \ \blambda \in \Lambda\}\,.
$$
\end{definition}

The following two observations motivate our interest in the function space $EP_{\Gamma,\Lambda}$.

\begin{remark} \label{rem:Q}

$(i)$ An exponential polynomial $p(\bx)=\bx^{\bgamma} e^{\blambda\, \bx} \in EP_{\Gamma,\Lambda}$  is a polynomial,
if $\blambda=\bnull$, or is an exponential function, if $\bgamma=\bnull$.  If $\bgamma=\bnull$ and $\blambda\in i\RR^s$,
then $p$ is a trigonometric function, or a hyperbolic function,
if $\bgamma=\bnull$ and $\blambda\in \RR^s$. The reproduction and generation of $E_{\Gamma,\Lambda}$ or its subspaces 
by non-stationary subdivision are important in CAGD, motion planning or iso-geometric analysis.

$(ii)$ We would like to emphasize that the definitions and the proofs of the results of this Section and Section
\ref{sec:algebra} still apply if one works with a subspace $EP_Q$ of $EP_{\Gamma,\Lambda}$,
where $Q \subset \Gamma \times \Lambda$ consists of pairs $(\bgamma,\blambda)$ for some $\bgamma \in \Gamma$ and $\blambda \in \Lambda$.
In this case, the algebraic conditions in \eqref{gen_zero_conditions}, \eqref{algebraic_conditions}, \eqref{algebraic_conditions_multiGENE},
\eqref{algebraic_conditions_multi} and \eqref{algebraic_conditions_multi_cor} should be checked for corresponding pairs 
$$
 (\bv,\bgamma), \quad \bv=\left( \epsilon_1 e^{-(\blambda \cdot M^{-(k+1)})_1}, \ldots, \epsilon_s e^{-(\blambda \cdot M^{-(k+1)})_s}
 \right), \quad \bepsilon \in \Xi, \quad (\blambda, \bgamma) \in Q.
$$ 
\end{remark}

\smallskip \noindent Since most of the properties of a subdivision scheme, e.g.
its convergence, smoothness or its support size, do
not depend on the choice of
$\bT^{[k]}=\{\bt_\balpha^{[k]}, \ \balpha \in \ZZ^s\}$, these are usually
set to
\begin{equation}\label{def:primal_par}
   \bt_\balpha^{[k]}=M^{-k} \balpha, \quad \balpha \in \ZZ^s, \quad k \ge 0.
\end{equation}
We refer to the choice in (\ref{def:primal_par}) as
\emph{standard} parametrization. On the contrary, the capability of subdivision to reproduce
exponential polynomials  does depend on the choice of
$\bt_\balpha^{[k]}$ and the standard parametrization is not always the optimal
one. We show in this section that the choice
\begin{equation}\label{def:general_par}
   \bt_\balpha^{[k]}=  M^{-k}(\balpha+ \btau), \quad \balpha \in \ZZ^s, \quad k\ge 0,
\end{equation}
with a suitable $\btau \in \RR^s$ turns out to be more advantageous.\\
We call the sequence $\{\bT^{[k]},\ k\ge 0\}$ with
$\bT^{[k]}=\{M^{-k}(\balpha+ \btau),\ \balpha\in \ZZ^s\}$  the \emph{parametrization
associated with a subdivision scheme} and $\btau \in \RR^s$  the corresponding \emph{shift} parameter.

\begin{definition} \label{def:ERgenerationlimit}
A convergent subdivision scheme $S_{\{\ba^{[k]},\ k\ge 0\}}$  is said to be
\emph{$EP_{\Gamma, \Lambda}$-generating}, if there exists a parametrization
$\{\bT^{[k]},\ k\ge 0\}$ with $\bT^{[k]}=\{\bt_{\balpha}^{[k]}=M^{-k}(\balpha+ \btau),
\ \balpha \in \ZZ^s\}$ and $\btau \in \RR^s$ such that
for every initial sequence $\bbf^{[0]}=\{ p(\bt_\balpha^{[0]}),\ \balpha \in \ZZ^s\}$,
$p \in EP_{\Gamma, \Lambda}$, we have
$$
 \lim_{k \rightarrow \infty} \, S_{\ba^{[\ell+k]}} S_{\ba^{[\ell+k-1]}} \ldots
  S_{\ba^{[\ell]}} \bbf^{[0]} \ \in \  EP_{\Gamma, \Lambda} \quad \forall \ell \geq 0.
$$
\end{definition}

We continue by defining  the notion of \emph{ $EP_{\Gamma, \Lambda}$-reproduction}.

\begin{definition}\label{def:ERreproductionlimit}
A convergent subdivision scheme $S_{\{\ba^{[k]},\ k\ge 0\}}$ is said to be
\emph{$EP_{\Gamma, \Lambda}$-reproducing}, if there exists a parametrization
$\{\bT^{[k]},\ k\ge 0\}$ with $\bT^{[k]}=\{\bt_{\balpha}^{[k]}=M^{-k}(\balpha+ \btau), \ \balpha \in \ZZ^s\}$ and $\btau \in \RR^s$ such that
for every initial sequence $\bbf^{[0]}=\{ p(\bt_\balpha^{[0]}), \ \balpha \in \ZZ^s\}$, $p \in EP_{\Gamma, \Lambda}$,
we have
$$
 \lim_{k \rightarrow \infty} \, S_{\ba^{[\ell+k]}} S_{\ba^{[\ell+k-1]}} \cdots
  S_{\ba^{[\ell]}} \bbf^{[0]} =p \quad \forall \ell \geq 0.
$$
\end{definition}

\noindent
Note that the generation and the reproduction properties are independent of the starting level $\ell$ of
refinement.

\smallskip \noindent We define next the  step-wise reproduction property of subdivision which  is
easier to check than its $EP_{\Gamma, \Lambda}$-reproduction.

\begin{definition} \label{def:ERreproduction_primal}
A convergent subdivision scheme $S_{\{\ba^{[k]},\ k\ge 0\}}$ is said to be
\emph{step-wise $EP_{\Gamma, \Lambda}$-reproducing}, if there exists a parametrization
$\{\bT^{[k]},\ k\ge 0\}$ with $\bT^{[k]}=\{\bt_{\balpha}^{[k]}=M^{-k}(\balpha+ \btau),\ \balpha \in \ZZ^s\}$
and $\btau \in \RR^s$ such that for the sequences
$\bbf^{[k]}=\{ p(\bt_\balpha^{[k]}), \ \balpha \in \ZZ^s\}$, $k \ge 0$, and $p \in EP_{\Gamma, \Lambda}$, we have
\begin{equation}\label{passo}
    \bbf^{[k+1]}=S_{\ba^{[k]}}\bbf^{[k]}.
\end{equation}
\end{definition}
It has been already observed in \cite{CC2012} in the stationary multivariate case and in \cite{CR2011} in
the binary non-stationary univariate setting
 that for non-singular schemes the concepts of  reproduction and step-wise reproduction are equivalent.
Since these results extend easily to the multivariate case we state the following
proposition without a proof.

\begin{proposition} \label{step-limit}
A non-singular non-stationary subdivision scheme 
is \emph{step-wise $EP_{\Gamma, \Lambda}$-reproducing}   if and only if it is
$EP_{\Gamma, \Lambda}$-reproducing.
\end{proposition}

\section{Algebraic conditions for generation and reproduction of exponential polynomials} \label{sec:algebra}

In this section we derive algebraic conditions  that guarantee
$EP_{\Gamma, \Lambda}-$reproduction by non-stationary subdivision schemes. These conditions are given in terms of
the symbols $\{a^{[k]}, \ k \ge 0\}$ which are evaluated at the elements of the sets
$$
 V_k=\{ (v_1, \ldots, v_s)^T \ : \ v_j= \epsilon_j e^{-(\blambda  \cdot M^{-(k+1)})_j},
  \  \blambda \in \Lambda, \ \bepsilon \in \Xi \ \}, \quad k \ge 0.
$$
For $k \ge 0$, we also define the following sets of vectors
$$
 V'_k=\{ (v_1, \ldots, v_s)^T \ : \ v_j= \epsilon_j e^{-(\blambda  \cdot M^{-(k+1)})_j}, \
    \blambda \in \Lambda, \ \bepsilon \in \Xi \setminus\{{\bf 1}\}\ \}.
$$

To simplify the presentation of the results presented in this section, in subsection \ref{subsec:uni_generation_reproduction},
we consider  the univariate case first
and then extend our results to the  multivariate setting in subsection \ref{subsec:multi_generation_reproduction}.

\subsection{Univariate case} \label{subsec:uni_generation_reproduction}

In the univariate case,  $M=m \ge 2$, $E=\{0,\ldots,m-1\}$ and
$\Xi =\{ e^{2\pi i m^{-1} \varepsilon } \ : \    \varepsilon \in E \}$ consists of the $m$-th
roots of unity. From \eqref{def:Soperator0}, we get that the $m$-ary subdivision scheme is  given by the
repeated application of
$m$ different rules
\begin{equation}\label{def:Soperator}
  \left(S_{\ba^{[k]}} \bbf^{[k]}\right)_{m\alpha+\varepsilon}=
   f^{[k+1]}_{m\alpha+\varepsilon} = \sum_{\beta \in \ZZ}
  \ra^{[k]}_{m\beta+\varepsilon} f^{[k]}_{\alpha-\beta},\quad \alpha \in\ZZ,
  \quad \varepsilon \in E,
 \quad k \ge 0 \,.
\end{equation}
 The structure of the sets $V_k,\ V'_k,\ k \ge 0, $ is a lot
simpler in the univariate case, namely,
$$
 V_k=\{ \epsilon\, e^{-\lambda m^{- (k+1)} } \ : \  \lambda \in \Lambda, \ \epsilon \in \Xi \},
$$
and
$$ V'_k=\{ \epsilon \, e^{-\lambda  m^{-(k+1)}} \ : \  \lambda \in \Lambda, \ \epsilon \in \Xi \setminus\{1\}\},
$$
for $k \ge 0$, respectively.

\begin{remark}
We remark that since our goal is to make the multivariate extension of
univariate results as straightforward as possible, in this subsection we use
notation less common for the univariate setting. For example, we write
$\displaystyle{\sum_{\varepsilon\in E}\rc_\varepsilon,\ E=\{0,\ldots,m-1\}}$
instead of   $\displaystyle{\sum_{j=0}^{m-1}\rc_j}$.
\end{remark}

\subsubsection{Generation of exponential polynomials}

The following result characterizes the $EP_{\Gamma, \Lambda}$-generation of a non-singular scheme in terms of
the so-called  zero conditions \eqref{gen_zero_conditions}. The proof of Proposition \ref{prop:EPgeneration} in the
case $m=2$ is given in \cite[Theorem 1]{VonBluUnser07}. We give the
generalization of this result to the case $m \ge 2$ in the notation familiar to subdivision audience. 

\begin{proposition}\label{prop:EPgeneration}
A  non-singular non-stationary subdivision scheme defined by the symbols $ \{ a^{[k]}(z),\  k\ge 0\}$
is
$EP_{\Gamma, \Lambda}$-generating if and only if
\begin{equation} \label{gen_zero_conditions}
 D^\gamma a^{[k]}(v)=0,\quad \gamma \in \Gamma, \quad v \in V'_k,\quad k\ge 0,
\end{equation}
for $V'_k=\{ \epsilon \, e^{-\lambda m^{-(k+1)}} \,: \,\lambda \in \Lambda, \,\epsilon \in \Xi \setminus\{1\}\}$.
\end{proposition}
\proof Let $\ell \ge 0$ and $\lambda \in \Lambda$.
We multiply both sides  of the non-stationary refinement equation \eqref{eq:ns-refinability}
by $e^{-\lambda m^{-\ell} x}$, $x \in \RR$, and get
\begin{equation}\label{eq:sopra}
    e^{-\lambda m^{-\ell} x}\phi^{[\ell]}(x)=\sum_{\alpha \in \ZZ}
    \ra^{[\ell]}_\alpha e^{-\lambda m^{-(\ell+1)}\alpha} \phi^{[\ell+1]}(mx-\alpha) e^{ -\lambda m^{-(\ell+1)} (mx-\alpha)}\,.
\end{equation}
Set  $\Phi_\ell(x)=e^{-\lambda m^{-\ell} x}\phi^{[\ell]}(x)$. Then, equation \eqref{eq:sopra} becomes
\begin{equation}
    \Phi_\ell(x)=\sum_{\alpha \in \ZZ} \ra^{[\ell]}_\alpha e^{-\lambda m^{-(\ell+1)}\alpha }\Phi_{\ell+1}(mx-\alpha)\,,
\end{equation}
or, equivalently, on the Fourier side,
\begin{equation}\label{eq:crucial}
    {\widehat \Phi}_\ell(\omega)=
     m^{-1} a^{[\ell]} \left(e^{ -\lambda m^{-(\ell+1)}-{i} m^{-1}\omega }\right)
     {\widehat \Phi}_{\ell+1}(m^{-1} \omega ),
\end{equation}
where $ \displaystyle a^{[\ell]}\left(e^{-{i}\omega} \right)
     =\sum_{\alpha \in \ZZ} \ra^{[\ell]}_\alpha e^{-{i} \omega \alpha}, \quad \omega \in \RR$.
The integer shifts of $\phi^{[\ell]}$ are $EP_{\Gamma, \Lambda}$-generating, i.e.
$$
 EP_{\Gamma, \Lambda} \subset \hbox{span} \{\phi^{[\ell]}(x-\alpha) \ : \  \alpha \in \ZZ\}\,,
$$
if and only if the functions $\Phi_\ell$ for all $\lambda \in \Lambda$ satisfy the so-called Strang-Fix conditions
$$
  D^\gamma {\widehat \Phi_\ell}(\omega) |_{\omega=2 \pi \beta}=0 \quad \hbox{for} \quad \beta \in E \setminus\{0\},
  \quad \gamma \in \Gamma,
$$
and ${\widehat \Phi}_\ell(0)\neq 0$, see e.g. \cite{StrangFix73}.

\smallskip
Therefore, evaluation of \eqref{eq:crucial} at $2\pi \beta$, $\beta \in E \setminus\{0\}$, yields
 \begin{equation}\label{eq:crucial2}
  0={\widehat \Phi}_{\ell}(2\pi \beta)  =m^{-1} a^{[\ell]} \left(
  e^{-\lambda m^{-(\ell+1)}-2 {i} \pi m^{-1} \beta }\right) {\widehat \Phi}_{\ell+1}(2 \pi m^{-1}\beta ).
\end{equation}
By Proposition \ref{prop:nonsingular_vs_independence}, the system of functions
$\{\Phi_\ell(\cdot-\alpha)\ :\ \alpha \in \ZZ\}$ is linearly independent for each $\ell \ge 0$. By
\cite{Ron89}, its linear independence is equivalent to the fact that the set $\{\xi \in \CC\ :\ {\widehat \Phi}_{\ell}(\xi+2\pi \beta)=0
\ \hbox{for all} \ \beta \in \ZZ\}$ is empty. Therefore, with $\xi=2 \pi m^{-1}\beta \in \CC$, the identities
\eqref{eq:crucial2} are satisfied if and only if
$$
 a^{[\ell]}(v)=0,\quad v \in V'_\ell,\quad \ell \ge 0\,.
$$
Taking  derivatives of both sides of \eqref{eq:crucial}
$$
D^\gamma \widehat{ \Phi_\ell} (\omega)=m^{-1} \sum_{\eta \le \gamma}  \left(
                                                                       \begin{array}{c}
                                                                         \gamma \\
                                                                         \eta \\
                                                                       \end{array}
                                                                     \right)
D^\eta a^{[\ell]}\left(e^{ -\lambda m^{-(\ell+1)}-im^{-1}\omega }\right)
  D^{\gamma-\eta} {\widehat \Phi_{\ell+1}}
(m^{-1}\omega)  
$$
and evaluating them at $2 \pi \beta$, $\beta \in E \setminus\{0\}$ yields, by induction on $\gamma$, that
the non-singular scheme is $EP_{\Gamma, \Lambda}$-generating if and only if
the identities in \eqref{gen_zero_conditions} are satisfied.
\qed

\medskip \noindent We would like to emphasize that the generation properties of subdivision schemes
are well-understood. Our
interest lies in better understanding of their reproduction properties.

\subsubsection{Reproduction of exponential polynomials} \label{subsec:reproduction_uni}

In this subsection, in Theorem \ref{th:characterization_poly_reproduction}, we derive algebraic conditions
on the mask symbols $\{a^{[k]}(z),\ k \ge 0\}$ that characterize the
$EP_{\Gamma, \Lambda}$-reproduction property of the associated non-stationary subdivision scheme.
We start by proving Proposition \ref{proposition:aux_step_wise_reproduction_uni}  that constitutes the main part
of the proof of Theorem \ref{th:characterization_poly_reproduction}.

\begin{proposition} \label{proposition:aux_step_wise_reproduction_uni}
A subdivision scheme $S_{\{ \ba^{[k]},\ k \ge 0\}}$  is step-wise
$EP_{\Gamma, \Lambda}$-reproducing if and only if there exists a shift parameter $\tau \in \RR$
such that
$$
  \sum_{\beta \in \ZZ} \ra^{[k]}_{m\beta+\varepsilon} \left(m \beta \right)^{\gamma'}  v^{m\beta+\varepsilon}=
   v^{m\tau-\tau} \epsilon^{-m\tau+\tau+\varepsilon} \left( m\tau-\tau-\varepsilon \right)^{\gamma'}
$$
is satisfied for all
$v \in V_k=\{ \epsilon \, e^{-\lambda m^{-(k+1)} } \ : \  \lambda \in \Lambda, \,
\epsilon \in \Xi \}$, $k\ge 0$, $\gamma' \le \gamma$, $\gamma \in \Gamma$, $\varepsilon
\in E$ and $\epsilon \in \Xi$.
\end{proposition}
\pf
 Let $\gamma \in \Gamma$, $\lambda \in \Lambda$ and define
$$
 f^{[k]}_{\alpha}= \left(m^{-k}(\alpha+\tau) \right)^{\gamma} \cdot e^{\lambda m^{-k}(\alpha+\tau)}, \quad \alpha \in \ZZ,
 \quad k \ge 0.
$$
 By \eqref{def:Soperator} and by definition \ref{def:ERreproduction_primal}, the step-wise reproduction of this sequence  is equivalent to the existence of $\tau \in \RR$ such that
 \begin{eqnarray}\label{aux1}
   f^{[k+1]}_{m\alpha+\varepsilon}=\sum_{\beta \in \ZZ} \ra^{[k]}_{m\beta+\varepsilon} f^{[k]}_{\alpha-\beta},
   \quad \varepsilon \in E, \quad \alpha \in \ZZ,\quad k\ge 0.
 \end{eqnarray}
 Multiplying both sides of \eqref{aux1}
 by $e^{-\lambda m^{-(k+1)}(m\alpha+m\tau+\varepsilon)}$ we get
 \begin{eqnarray} \label{aux1_1}
   &&\hspace{-2cm} \sum_{\beta \in \ZZ} \ra^{[k]}_{m\beta+\varepsilon} \left(m^{-k}(\alpha-\beta+ \tau) \right)^{\gamma}
   e^{-\lambda m^{-(k+1)}( m \beta+\varepsilon)} \notag
   \\ && = \left(m^{-(k+1)}(m\alpha+\varepsilon +\tau) \right)^{\gamma}  e^{-\lambda m^{-(k+1)}( m \tau-\tau)},
 \end{eqnarray}
or, equivalently,
 \begin{eqnarray*}
   &&\hspace{-2cm} \sum_{\beta \in \ZZ} \ra^{[k]}_{m\beta+\varepsilon} \left(m^{-k}(\alpha-\beta+ \tau) \right)^{\gamma}
   e^{(-\lambda m^{-(k+1)}-2\pi {i} m^{-1} \varepsilon)
\cdot ( m \beta+\varepsilon)} \notag
   \\ && =e^{-2\pi  i m^{-1} \varepsilon \cdot (m \beta+\varepsilon-m\tau+\tau)}  \left(m^{-(k+1)}(m\alpha+\varepsilon
    +\tau) \right)^{\gamma}  e^{(-\lambda m^{-(k+1)}-2\pi {i} m^{-1} \varepsilon)( m \tau-\tau)}.
 \end{eqnarray*}
Due to the fact that
 $$
  e^{2\pi {i}  m^{-1} \varepsilon \cdot m\beta}=1 \quad \hbox{for all} \quad \varepsilon \
  \in E \quad \hbox{and} \quad \beta \in \ZZ,
 $$
we have that
 \eqref{aux1_1} is equivalent to
 \begin{eqnarray} \label{aux2}
   &&\hspace{-2cm} \sum_{\beta \in \ZZ} \ra^{[k]}_{m\beta+\varepsilon} \left(m^{-k}(\alpha-\beta+ \tau) \right)^{\gamma}
   v^{ m \beta+\varepsilon} \notag
   \\ && =
    \epsilon^{- m \tau+\tau+\varepsilon} \left(m^{-(k+1)}(m\alpha+\varepsilon +\tau) \right)^{\gamma}  v^{ m \tau-\tau},
 \end{eqnarray}
 for $v \in V_k$. Considering both sides of \eqref{aux2} as polynomials in $m^{-k}\alpha$ with
 $$
    \left(m^{-k}(\alpha-\beta+ \tau) \right)^{\gamma}=\sum_{\eta \le \gamma} \left(\begin{array}{c} \gamma \\ \eta \end{array} \right)
    \left(m^{-k}(-\beta+ \tau) \right)^{\eta} \left(m^{-k} \alpha \right)^{\gamma-\eta}\,,
 $$
 and
 $$
    \left(m^{-(k+1)}(m\alpha+\varepsilon+ \tau) \right)^{\gamma}=\sum_{\eta \le \gamma} \left(\begin{array}{c} \gamma \\ \eta
     \end{array} \right)
    \left(m^{-(k+1)}(\varepsilon+ \tau) \right)^{\eta} \left(m^{-k} \alpha \right)^{\gamma-\eta},
 $$
 we get that \eqref{aux2} is equivalent to
  \begin{eqnarray} \label{aux3}
   && \hspace{-1cm}\sum_{\beta \in \ZZ} \ra^{[k]}_{m\beta+\varepsilon} \left(m^{-k}(-\beta+ \tau) \right)^{\eta}
   v^{ m \beta+\varepsilon} \notag
   \\ &&=
    \epsilon^{ -m\tau+\tau+\varepsilon} \left(m^{-(k+1)}(\varepsilon +\tau) \right)^{\eta}  v^{ m \tau-\tau}, \quad
    v \in V_k, \quad \eta \le \gamma.
 \end{eqnarray}
 Next, we expand
 $$
   \left(m^{-(k+1)}(m\tau+ \tau-m\tau+\varepsilon) \right)^{\eta}=\sum_{\gamma' \le \eta} \left(\begin{array}{c}
   \eta \\ \gamma' \end{array} \right)
    \left(m^{-(k+1)}(\tau-m\tau+\varepsilon) \right)^{\gamma'} \left(m^{-k} \tau \right)^{\eta-\gamma'}\,,
 $$
 and, similarly, $ \left(m^{-k}(-\beta+ \tau) \right)^{\eta}$. Thus, \eqref{aux3} is satisfied if and only if
  \begin{eqnarray} \label{aux4}
   && \hspace{-1cm}\sum_{\beta \in \ZZ} \ra^{[k]}_{m\beta+\varepsilon} \left(-m^{-k}\beta \right)^{\gamma'} v^{ m \beta+\varepsilon}
   =\notag
   \\ && \hspace{-0.5cm}
    \epsilon^{- m \tau+\tau+\varepsilon} \left(m^{-(k+1)}(\tau +\varepsilon -m\tau) \right)^{\gamma'}  v^{ m\tau-\tau},
    \
    v \in V_k, \  \gamma' \le \eta \le \gamma.
 \end{eqnarray}
 Multiplying both sides by $(-m^{k+1})^{\gamma'}$, we get the claim.
\eop

\noindent To formulate the main result of this subsection we  define
\begin{equation} \label{def:q}
q_{0}(z)= 1,\quad  q_{\gamma}(z)=  \prod_{j=0}^{\gamma-1} (z-j), \quad \gamma\in \NN_0, \quad z \in \CC.
\end{equation}

\begin{theorem}\label{th:characterization_poly_reproduction}
A  non-singular subdivision scheme   $S_{\{\ba^{[k]}, \  k \ge 0\}}$
reproduces $EP_{\Gamma, \Lambda}$ if and only if there exists a shift parameter $\tau \in \RR$ such that
\begin{equation} \label{algebraic_conditions}
 v^\gamma D^\gamma a^{[k]}(v)= \left\{
 \begin{array}{ll}
 m \cdot v^{m\tau-\tau}  q_\gamma(m\tau-\tau),& \hbox{for all}\ v\  \hbox{such that}\  \epsilon=1, \\
 0,& \hbox{otherwise}, \end{array} \right.
\end{equation}
for  all $v \in V_k=\{ \epsilon\, e^{-\lambda m^{-(k+1)} } \ : \   \lambda \in \Lambda, \  \epsilon \in \Xi \}
,\ k\ge 0$, $\gamma \in \Gamma$.
\end{theorem}
\pf
Due to Proposition \ref{step-limit} it suffices to prove this statement for step-wise $EP_{\Gamma, \Lambda}$-reproduction.
We use Proposition \ref{proposition:aux_step_wise_reproduction_uni} and show that the step-wise $EP_{\Gamma, \Lambda}$-reproduction,
is equivalent to the identities \eqref{algebraic_conditions}.
To this purpose we consider
$q_\gamma(m\alpha+\varepsilon)$, for fixed $\gamma \in \Gamma$ and
$\varepsilon \in E$, as a polynomial in $ m \alpha$
and write
\begin{equation}\label{aux5}
   q_\gamma(m \alpha+\varepsilon)=\sum_{\eta \le \gamma} c_{\gamma,\varepsilon,\eta} (m \alpha)^\eta,
   \quad c_{\gamma,\varepsilon,\eta} \in \RR.
\end{equation}
From \eqref{eq:decomp} we get
$$
 z^\gamma D^\gamma a^{[k]}(z)=\sum_{\varepsilon \in E} \sum_{\alpha \in \ZZ} \ra^{[k]}_{m\alpha+\varepsilon} \
 q_\gamma(m\alpha+\varepsilon) z^{m\alpha+\varepsilon}, \quad \gamma \in \Gamma.
$$
Let $v \in V_k$. By Proposition \ref{proposition:aux_step_wise_reproduction_uni} and by \eqref{aux5}, the step-wise $EP_{\Gamma, \Lambda}$-reproduction is equivalent to
\begin{eqnarray*}
 v^\gamma D^\gamma a^{[k]}(v) &=&\sum_{\varepsilon \in E} \sum_{\alpha \in \ZZ} \ra^{[k]}_{m\alpha+\varepsilon} \
 q_{\gamma}(m\alpha+\varepsilon) \ v^{m\alpha+\varepsilon} \\
 &=& \sum_{\varepsilon \in E} \sum_{\eta \le \gamma} c_{\gamma,\varepsilon,\eta} \sum_{\alpha \in \ZZ} \ra^{[k]}_{m\alpha+
 \varepsilon} (m\alpha)^\eta v^{m\alpha+\varepsilon}\\
 &=& v^{m\tau-\tau} \sum_{\varepsilon \in E}  \epsilon^{-m\tau+\tau+\varepsilon} \sum_{\eta \le \gamma} c_{\gamma,\varepsilon,\eta}  \left( m\tau -\tau-\varepsilon
 \right)^\eta \\
 &=& v^{m\tau-\tau} \epsilon^{-m\tau+\tau} q_\gamma(m\tau-\tau) \sum_{\varepsilon \in E} \epsilon^{\varepsilon}
\end{eqnarray*}
for some $\tau \in \RR$. The claim follows, due to
$$
 \sum_{\varepsilon \in E} \epsilon^{\varepsilon}=\left\{
                                  \begin{array}{ll}
                                    m, & \epsilon=1, \\
                                     0, & \hbox{otherwise.}
                                   \end{array}
                                 \right. $$
\eop



\subsection{Multivariate case} \label{subsec:multi_generation_reproduction}

In this subsection we give a closer look at the multivariate case. We will not repeat
results that can be easily extended from the univariate setting by simply
replacing $\ZZ,\ \RR$ or $\CC$ with $\ZZ^s,\ \RR^s$ or $\CC^s$, respectively. This is the case of
Proposition \ref{proposition:generation_multi} which we give without proof.

\smallskip \noindent
Recall that $m=|det(M)|$ determines the cardinality of
$E=\{\bvarepsilon_0, \ldots, \bvarepsilon_{m-1}\}$ and of  $\Xi =\{ e^{2\pi  i  M^{-T}\bxi} \ : \   \bxi
\ \hbox{is a coset representative of} \ M^{-T}\ZZ^s / \ZZ^s \}$.

\begin{proposition} \label{proposition:generation_multi}  A subdivision scheme  defined by
the symbols
$\{a^{[k]}(\bz), \ k \ge 0\}$
is $EP_{\Gamma, \Lambda}$-generating if and only if
\begin{equation} \label{algebraic_conditions_multiGENE}
D^\bgamma a^{[k]}(\bv)=0,\ \quad \bv \in V'_k,\quad \bgamma \in \Gamma, \quad k\ge 0\,,
\end{equation}
for $
 V'_k=\{ (v_1,\ldots,v_s)^T \ : \ v_j=\epsilon_j e^{-(\blambda  \cdot M^{-(k+1)})_j},
  \  \blambda \in \Lambda, \ \ \bepsilon \in \Xi \setminus\{{\bf 1}\}\ \}\,.
$
\end{proposition}

\subsubsection{Reproduction of multivariate exponential polynomials } \label{subsec:reproduction_multi}

The  multivariate extension of Proposition
\ref{proposition:aux_step_wise_reproduction_uni} can appear trivial after
reading the proof of Proposition
\ref{proposition:aux_step_wise_reproduction_multi}. We believe that would not be
the case, if we omitted its proof. There are also several crucial differences
between the proofs of Theorem \ref{th:characterization_poly_reproduction_multi} and
its univariate counterpart.

\begin{proposition} \label{proposition:aux_step_wise_reproduction_multi}  A subdivision scheme
$S_{\{\ba^{[k]}, \ k \ge 0\}}$ is step-wise $EP_{\Gamma, \Lambda}$-reproducing
if and only if there exists a shift parameter $\btau \in \RR^s$ such that
$$
  \sum_{\bbeta \in \ZZ^s} \ra^{[k]}_{M\bbeta+\bvarepsilon} \left(M^{-k}\bbeta \right)^{\bgamma'}
   \bv^{M\bbeta+\bvarepsilon}=
   \bv^{M\btau-\btau} \bepsilon^{-M\btau+\btau+\bvarepsilon} \left(M^{-(k+1)}(M\btau-\btau-\bvarepsilon) \right)^{\bgamma'}
$$
is satisfied for all $\bv \in V_k=\{ (v_1,\ldots,v_s)^T\ : \ v_j=\epsilon_j e^{-(\blambda  \cdot M^{-(k+1)})_j}, \ \blambda \in \Lambda,\,\bepsilon \in \Xi \}$, $k\ge 0$, $\bgamma' \le \bgamma$, $\bgamma \in \Gamma$, $\bvarepsilon \in E$ and $\bepsilon \in \Xi$.
\end{proposition}
\pf
Let $\bgamma \in \Gamma$, $\blambda \in \Lambda$ and define
$$
 f^{[k]}_{\balpha}= \left(M^{-k}(\balpha+\btau) \right)^{\bgamma} \cdot e^{\blambda \cdot M^{-k}(\balpha+\btau)},
  \quad \balpha \in \ZZ^s, \quad k \ge 0.
$$
 By definition, the step-wise reproduction of sequences sampled from the exponential-polynomials in $EP_{\Gamma, \Lambda}$ is equivalent to the existence of $\btau \in \RR^s$ such that
 \begin{eqnarray}\label{aux1_m}
   f^{[k+1]}_{M\balpha+\bvarepsilon}=\sum_{\bbeta \in \ZZ^s} \ra^{[k]}_{M\bbeta+\bvarepsilon}
   f^{[k]}_{\balpha-\bbeta},  \quad \bvarepsilon \in E, \quad \balpha \in \ZZ^s.
 \end{eqnarray}
 Multiplying both sides of \eqref{aux1_m}
 by $e^{-\blambda \cdot M^{-(k+1)}(M\balpha+M\btau+\bvarepsilon)}$ we get
 \begin{eqnarray} \label{aux1_1_m}
   &&\hspace{-2cm} \sum_{\bbeta \in \ZZ^s} \ra^{[k]}_{M\bbeta+\bvarepsilon}
   \left(M^{-k}(\balpha-\bbeta+ \btau) \right)^{\bgamma}  e^{-\blambda \cdot M^{-(k+1)}( M \bbeta+\bvarepsilon)} \notag
   \\ && = \left(M^{-(k+1)}(M\balpha+\bvarepsilon +\btau) \right)^{\bgamma}  e^{-\blambda \cdot M^{-(k+1)}( M \btau-\btau)},
 \end{eqnarray}
or, equivalently,
 \begin{eqnarray*}
   &&\hspace{-2cm} \sum_{\bbeta \in \ZZ^s} \ra^{[k]}_{M\bbeta+\bvarepsilon}
   \left(M^{-k}(\balpha-\bbeta+ \btau) \right)^{\bgamma}
   e^{(-\blambda \cdot M^{-(k+1)}-2\pi   i M^{-T} \bxi)
\cdot ( M \bbeta+\bvarepsilon)} \notag
   \\ && =e^{-2\pi  i M^{-T} \bxi \cdot (M \bbeta+\bvarepsilon-M\btau+\btau)}
   \left(M^{-(k+1)}(M\balpha+\bvarepsilon +\btau) \right)^{\bgamma}
   e^{(-\blambda \cdot M^{-(k+1)}-2\pi {i} M^{-T} \bxi) \cdot ( M \btau-\btau)}.
 \end{eqnarray*}
Now, by properties of scalar products we have $ M^{-T} \bxi \cdot M \bbeta=
\bxi \cdot M^{-1} M \bbeta=\bxi \cdot \bbeta \in \ZZ$, and therefore
$e^{-2\pi i M^{-T} \bxi \cdot M \bbeta}=1$.
Thus,
\eqref{aux1_1_m} is equivalent to
 \begin{eqnarray} \label{aux2_m}
   &&\hspace{-2cm} \sum_{\bbeta \in \ZZ^s} \ra^{[k]}_{M\bbeta+\bvarepsilon} \left(M^{-k}(\balpha-\bbeta+ \btau)
   \right)^{\bgamma} \bv^{ M \bbeta+\bvarepsilon} \notag
   \\ && =
    \bepsilon^{- M \btau+\btau+\bvarepsilon} \left(M^{-(k+1)}(M\balpha+\bvarepsilon +\btau) \right)^{\bgamma}  \bv^{ M \btau-\btau},
 \end{eqnarray}
 for $\bv \in V_k$. Considering both sides of \eqref{aux2_m} as polynomials in $M^{-k}\balpha$ with
 $$
    \left(M^{-k}(\balpha-\bbeta+ \btau) \right)^{\bgamma}=\sum_{\bfeta \le \bgamma} \left(\begin{array}{c} \bgamma \\ \bfeta \end{array} \right)
    \left(M^{-k}(-\bbeta+ \btau) \right)^{\bfeta} \left(M^{-k} \balpha \right)^{\bgamma-\bfeta}
 $$
 and
 $$
    \left(M^{-(k+1)}(M\balpha+\bvarepsilon+ \btau) \right)^{\bgamma}=\sum_{\bfeta \le \bgamma} \left(\begin{array}{c}
    \bgamma \\ \bfeta \end{array} \right)
    \left(M^{-(k+1)}(\bvarepsilon+ \btau) \right)^{\bfeta} \left(M^{-k} \balpha \right)^{\bgamma-\bfeta},
 $$
 we get that \eqref{aux2_m} is equivalent to
  \begin{eqnarray} \label{aux3_m}
   && \hspace{-1cm}\sum_{\bbeta \in \ZZ^s} \ra^{[k]}_{M\bbeta+\bvarepsilon} \left(M^{-k}(-\bbeta+ \btau) \right)^{\bfeta}  \bv^{ M \bbeta+\bvarepsilon} \notag
   \\ &&=
    \bepsilon^{ -M \btau+\btau+\bvarepsilon} \left(M^{-(k+1)}(\bvarepsilon +\btau) \right)^{\bfeta}  \bv^{ M \btau-\btau}, \quad
    \bv \in V_k, \quad \bfeta \le \bgamma.
 \end{eqnarray}
 Next, we expand
 $$
   \left(M^{-(k+1)}(M\btau+ \btau-M\btau+\bvarepsilon) \right)^{\bfeta}=\sum_{\bgamma' \le \bfeta} \left(\begin{array}{c} \bfeta \\ \bgamma' \end{array} \right)
    \left(M^{-(k+1)}(\btau-M\btau+\bvarepsilon) \right)^{\bgamma'} \left(M^{-k} \btau \right)^{\bfeta-\bgamma'}
 $$
 and, similarly, $ \left(M^{-k}(-\bbeta+ \btau) \right)^{\bfeta}$. Thus, \eqref{aux3_m} is satisfied if and only if
  \begin{eqnarray*} \label{aux4_m}
   && \hspace{-1cm}\sum_{\bbeta \in \ZZ^s} \ra^{[k]}_{M\bbeta+\bvarepsilon} \left(- M^{-k}\bbeta \right)^{\bgamma'}
   \bv^{ M \bbeta+\bvarepsilon} =\notag
   \\ && \hspace{-0.5cm}
    \bepsilon^{- M \btau+\btau+\bvarepsilon} \left(M^{-(k+1)}(\btau +\bvarepsilon -M \btau) \right)^{\bgamma'}  \bv^{ M \btau-\btau}, \quad
    \bv \in V_k, \quad  \bgamma' \le \bfeta \le \bgamma.
 \end{eqnarray*}
 Multiplying both sides by $(-1)^{|\bgamma'|}$, we get the claim.
\eop

\smallskip \noindent
Define
\begin{equation} \label{def:q_m}
  q_{\bf 0}(z_1,\ldots,z_s)=1,\quad q_{\bgamma}(z_1,\ldots,z_s)=\prod_{\ell=1}^s \prod_{j=0}^{\gamma_\ell-1} (z_\ell-j), \quad \bgamma \in \NN_0^s.
\end{equation}

\begin{theorem}\label{th:characterization_poly_reproduction_multi}
A non-singular subdivision scheme $S_{\{\ba^{[k]}, \ k \ge 0\}}$  reproduces $EP_{\Gamma, \Lambda}$
if and only if there exists a shift parameter  $\btau \in \RR^s$  such that
\begin{equation} \label{algebraic_conditions_multi}
 \bv^\bgamma D^\bgamma a^{[k]}(\bv)= \left\{
 \begin{array}{ll}
 m \cdot \bv^{M\tau-\tau}  q_\bgamma(M\btau-\btau), &\hbox{for all}\ \bv\ \hbox{such that}\ \bepsilon=\bone, \\
 0,& \hbox{otherwise}, \end{array} \right.
\end{equation}
for  all $\bv \in V_k=\{  (v_1,\ldots,v_s)^T\ : \ v_j=\epsilon_j e^{-(\blambda  \cdot M^{-(k+1)})_j}, \ \blambda \in \Lambda, \, \bepsilon \in \Xi\}$,
$k\ge 0$, $\bgamma \in \Gamma$.
\end{theorem}
\pf
Due to Proposition \ref{step-limit} it suffices to prove this statement for step-wise reproduction of sequences sampled from elements of $EP_{\Gamma, \Lambda}$.
We use Proposition \ref{proposition:aux_step_wise_reproduction_multi} and show that the step-wise
$EP_{\Gamma, \Lambda}$-reproduction
is equivalent to the identities \eqref{algebraic_conditions_multi}.
Consider
$q_\bgamma(M\balpha+\bvarepsilon)=q_\bgamma(M^{k+1} M^{-k}\balpha+\bvarepsilon)$, for fixed $\bgamma \in \Gamma$
 and $\bvarepsilon \in E$, as a polynomial in $M^{-k}\balpha$, $\balpha \in \ZZ^s$,
and write
\begin{equation}\label{aux5_m}
   q_\bgamma(M^{k+1} M^{-k}\balpha+\bvarepsilon)=\sum_{\bfeta \le \bgamma} c_{\bgamma,\bvarepsilon,\bfeta}
   \left( M^{-k}\balpha \right)^\bfeta, \quad c_{\bgamma,\bvarepsilon,\bfeta} \in \RR.
\end{equation}
Let $\bv \in V_k$. By Proposition \ref{proposition:aux_step_wise_reproduction_multi} and by \eqref{aux5_m}, the step-wise reproduction is equivalent to
\begin{eqnarray*}
 \bv^\bgamma D^\bgamma a^{[k]}(\bv) &=&\sum_{\bvarepsilon \in E} \sum_{\balpha \in \ZZ^s} \ra^{[k]}_{M\balpha+\bvarepsilon} \
 q_{\bgamma}(M\balpha+\bvarepsilon) \bv^{M\balpha+\bvarepsilon} \\
 &=& \sum_{\bvarepsilon \in E} \sum_{\bfeta \le \bgamma} c_{\bgamma,\bvarepsilon,\bfeta} \sum_{\balpha \in \ZZ^s}
 \ra^{[k]}_{M\balpha+\bvarepsilon} \left( M^{-k}\balpha \right)^\bfeta \bv^{M\balpha+\bvarepsilon}\\
 &=& \sum_{\bvarepsilon \in E} \sum_{\bfeta \le \bgamma} c_{\bgamma,\bvarepsilon,\bfeta}
 \left( M^{-(k+1)}(M\btau -\btau-\bvarepsilon) \right)^\bfeta \bv^{M\btau-\btau} \bepsilon^{-M\btau+\btau+\bvarepsilon}\\
 &=& \bv^{M\btau-\btau} \bepsilon^{-M\btau+\btau} q_\bgamma(M\btau-\btau) \sum_{\bvarepsilon \in E} \bepsilon^{\bvarepsilon}
\end{eqnarray*}
for some $\btau \in \RR^s$. The claim follows due to
$$
 \sum_{\bvarepsilon \in E} \bepsilon^{\bvarepsilon}=\left\{
                                   \begin{array}{ll}
                                     m, & \bepsilon=\bone, \\
                                     0, & \hbox{otherwise.}
                                   \end{array}
                                 \right.
$$
\eop

In the case the scheme
$S_{\{\ba^{[k]},\ k \ge 0\}}$ is convergent, but not non-singular,
the conditions \eqref{algebraic_conditions_multi}
in Theorem \ref{th:characterization_poly_reproduction_multi}
are only sufficient for the $EP_{\Gamma, \Lambda}$-reproduction. Subdivision schemes
generating exponential box splines, whose translates are linearly dependent,
are important examples of such convergent schemes. See
\cite[Theorem 4.3]{Ron92} for conditions that characterize the
linear independence of exponential box splines.

\begin{corollary}  \label{cor:multi_repr}
A convergent subdivision scheme  $S_{\{\ba^{[k]}, \ k \ge 0\}}$ reproduces $EP_{\Gamma, \Lambda}$
if there exists a shift parameter $\btau \in \RR^s$ such that
\begin{equation} \label{algebraic_conditions_multi_cor}
 \bv^\bgamma D^\bgamma a^{[k]}(\bv)= \left\{
 \begin{array}{ll}
 m \cdot \bv^{M\tau-\tau}  q_\bgamma(M\btau-\btau), &\hbox{for all}\ \bv\ \hbox{such that}\ \bepsilon=\bone, \\
 0,& \hbox{otherwise}, \end{array} \right.
\end{equation}
for all $\bv \in \{ (v_1,\ldots,v_s)^T\ : \ v_j=\epsilon_j e^{-(\blambda  \cdot M^{-(k+1)})_j},  \, \blambda \in \Lambda, \, \bepsilon \in \Xi\}$,
$k\ge 0$, $\bgamma \in \Gamma$.
\end{corollary}

\section{Examples and applications} \label{sec:examples}

\subsection{Shift factors and interpolatory schemes}

The first important application of
Theorem \ref{th:characterization_poly_reproduction_multi}, as in the binary univariate non-stationary case or as in the
stationary multivariate case, is
the analysis of the reproduction properties of interpolatory schemes, i.e. the
schemes whose masks satisfy
$$
 \ra^{[k]}_{M\balpha}=\delta_{\balpha,0},\ \balpha\in \ZZ^s.
$$
It is also of importance to analyse the effect of the mask shifts on the reproduction properties
of the corresponding schemes. The following results
are easily obtained  by combining \cite[Corollary 2]{CR2011} and
\cite[Proposition 3.4]{CC2012} or \cite[Proposition
3]{CR2011} and \cite[Lemma 3.1]{CC2012}, respectively.

\begin{corollary}\label{cor:interp}
A non-singular interpolatory
subdivision scheme $S_{\{\ba^{[k]},\ k \ge 0\}}$ is $EP_{\Gamma, \Lambda}$-reproducing
only if  $\btau=\bnull$. Moreover, the generation of $EP_{\Gamma, \Lambda}$
implies that the scheme reproduces the same space $EP_{\Gamma, \Lambda}$.
\end{corollary}

\begin{corollary}\label{cor:shift}
If a non-singular subdivision scheme  $S_{\{\ba^{(k)},\ k\ge 0\}}$
reproduces  $EP_{\Gamma, \Lambda}$, then so does the
scheme $S_{\{\bb^{(k)},\ k\ge 0\}}$ defined by the symbols
$b^{(k)}(\bz)=\bz^\bbeta \cdot a^{(k)}(\bz),\ k\ge 0$ and $\bbeta \in\ZZ^s$.
\end{corollary}

\subsection{Subdivision schemes for exponential B-splines  and  exponential box splines} \label{subsec:expBsplines}

In this subsection we study the reproduction properties of the subdivision
schemes associated with the so-called exponential B-splines and box splines. The
symbols of these schemes even in the non-stationary case play a role of
smoothing factors \cite{DL95}. They also determine the generation properties of
subdivision schemes.

\subsubsection{Exponential B-splines}

It is well known that in the binary case   a non-stationary subdivision scheme generates univariate
exponential polynomials
\begin{equation}\label{exp-pol}
    p(x)=x^\gamma e^{\lambda x} \qquad \lambda\in \CC,\quad  \gamma \in \NN_0\,,
\end{equation}
if its symbol contains factors of the type
$$
  (1+r_k\,z) \quad \hbox{with a suitable} \  r_k\in \CC.
$$
We show that in the $m-$ary case such exponential polynomials $p$ are generated by non-stationary
subdivision schemes with symbols containing factors of the type
$$
 (1+r_k\,z+r^2_k\,z^2+\cdots+r^{m-1}_k\,z^{m-1}) \quad \hbox{with a suitable} \  r_k\in \CC.
$$
\smallskip \noindent The following result is a generalization of a result in \cite{DL95}.
We present it here as it also allows for derivation of masks of exponential box
splines in the next subsection and  illustrates the generation property of the corresponding schemes.

\begin{proposition} \label{prop:ex_bsplines}
Let a non-stationary $m-$ary subdivision scheme be given by
\begin{equation}\label{symbol:exp-poluni}
  a^{[k]}(z)=(1+r_k\,z+r^2_k\,z^2+\cdots+r^{m-1}_k\,z^{m-1}),\  r_k=e^{\lambda m^{-(k+1)}}, \ \lambda \in \CC,
  \ \ k\ge 0\,,
\end{equation}
then its basic limit function  is  $\phi(x)=e^{\lambda x} \chi_{[0,1)}$.
\end{proposition}

\pf
\smallskip \noindent First we observe that a subdivision scheme based on the masks \eqref{symbol:exp-poluni} with
 support size $N_k=m$, whenever convergent or weakly convergent, generates a basic limit function supported
 on $[0,\frac{N_k-1}{m-1}]=[0,1]$ , see for example \cite{ContiHormann11}.  The subdivision rules corresponding
 to the symbols in \eqref{symbol:exp-poluni}  are given by
 \begin{equation}\label{rules:exp-pol}
                     f^{k+1}_{m\alpha+\varepsilon}=r_k^\varepsilon f^k_\alpha,\quad
                     \varepsilon \in E \,.
 \end{equation}
 Starting the subdivision process with $\bdelta$, from \eqref{rules:exp-pol} we get that the value
 of the basic limit function $\phi$ at any $m$-adic point
 $$
     x=\sum_{j=1}^k m^{-j} \varepsilon_j,\quad \varepsilon_j \in E \,,
$$
is
$$
 \phi(x)=\prod_{j=1}^k r_{j-1}^{\varepsilon_j}\,.
$$
At every  non $m$-adic point $ \displaystyle x=\sum_{j=1}^\infty m^{-j} \varepsilon_j$,
$\varepsilon_j \in E$, we can define
$$
\phi(x)=\prod_{j=1}^\infty r_{j-1}^{\varepsilon_j},
$$
due to $\displaystyle \sum_{k \in \ZZ_+} |1-r_k| < \infty$, which ensures the convergence of this infinite product. 
For the point $x \not \in [0,1)$ we
set $\phi(x)=0$. To show the continuity of $\phi$ at non $m$-adic points we consider a sequence of points
$$
 \{x^\ell,\ \ell \in \NN_0 \}, \quad \hbox{with}\quad x^\ell=\sum_{j=1}^{k_\ell} m^{-j} \varepsilon_j
  +\sum_{j=k_\ell+1}^\infty m^{-j} \varepsilon^\ell_j, \quad \varepsilon^\ell_j \in E,
$$
such that $k_\ell$ goes to infinity as $\ell$ goes to infinity. Then
$$
\lim_{\ell \rightarrow\infty}|\phi(x^\ell)-\phi(x)|=\lim_{\ell \rightarrow\infty}\left|\prod_{j=1}^\infty r_{j-1}^{\varepsilon_j}-
 \prod_{j=1}^{k_\ell} r_{j-1}^{\varepsilon_j} \prod_{j=k_\ell+1}^\infty r_{j-1}^{\varepsilon^\ell_j} \right| =0.
$$
The non-uniqueness of the representations of $m$-adic points
$$
   x= \sum_{j=1}^k m^{-j} \varepsilon_j = \sum_{j=1}^{k-1} m^{-j} \varepsilon_j +
    m^{-k} (\varepsilon_k-1) +  \sum_{j=k+1}^{\infty} m^{-j} (m-1)
$$
with
$$
 \sum_{j=k+1}^\infty m^{-j} (m-1)= m^{-k},
$$
makes the analysis of continuity at these points more involved.
Thus, we need to consider two types of sequences  that converge to $x$ from the left and from the right. This
implies that $\phi$ is continuous at $m$-adic points   if and only if
$$
   r_{0}^{\varepsilon_1} \ldots r_{k-2}^{\varepsilon_{k-1}} \left( r_{k-1}^{\varepsilon_{k}} -
   r_{k-1}^{\varepsilon_{k}-1} \prod_{k+1}^\infty r_{j-1}^{m-1} \right)=0,
$$
or, equivalently,
$$
r_{k-1}=\prod_{j=k+1}^\infty r_{j-1}^{m-1}=r_{k}^{m-1}\prod_{j=k+2}^\infty r_{j-1}^{m-1}=
r_{k}^m\quad \Leftrightarrow \quad r_k=r_{k-1}^{1/m}=\cdots=r_{0}^{1/m^k}\,.
$$
Hence, for $r_0=e^{\lambda/m}$, $\lambda \in \CC$, and $r_k=e^{\lambda/m^{k+1}}$ we get
$$
\phi(x)=\prod_{j=1}^\infty r_{j-1}^{\varepsilon_j}=\prod_{j=1}^\infty r_{0}^{\varepsilon_j m^{-(j+1)}}
=\prod_{j=1}^\infty \left(e^{ \lambda \varepsilon_j m^{-j}}\right)
=e^{\lambda    \sum_{j=1}^\infty  m^{-j} \varepsilon_j}=e^{\lambda x}.
$$
\eop

Proposition \ref{prop:ex_bsplines} implies that the corresponding subdivision scheme generates
the exponential polynomials $e^{\lambda x}$. Next, we use the algebraic conditions on
the subdivision symbols in Theorem \ref{th:characterization_poly_reproduction}  to show how
to normalize the symbols appropriately to ensure the reproduction of exponential
polynomials $e^{\lambda x}$.

\begin{lemma}
Let $\Gamma=\{0\}$ and $\Lambda=\{\lambda\}$, $\lambda \in \CC$.
Then a non-singular non-stationary scheme given by
\begin{equation}\label{eq:simbol1f}
    a^{[k]}(z)=K^{[k]} (1+r_k z+r^2_k z^2+\cdots+r^{m-1}_k z^{m-1}), \ r_k=e^{\lambda
 m^{-(k+1)}},\ k \ge 0\,,
\end{equation}
reproduces $EP_{\Gamma, \Lambda}=\{e^{\lambda x}\}$ if and only if $K^{[k]}=r_k^{-m\tau+\tau}$ and $\tau \in \RR$.
\end{lemma}
\pf
In this case
$$
 V_k=\{ e^{2\pi\, i\varepsilon/m} e^{-\lambda\,m^{-(k+1)}}\ : \ \varepsilon\in E\},\quad  k\ge 0\,.
$$
We only need to check the algebraic condition $a^{[k]}(v)=mv^{m\tau-\tau}$ for $v=r_k^{-1}=e^{-\lambda\,m^{-(k+1)}}$ as
the rest of the conditions in \eqref{algebraic_conditions} are trivially satisfied.
By Theorem \ref{th:characterization_poly_reproduction} the scheme reproduces $EP_{\Gamma, \Lambda}$
if and only if for all $k\ge 0$ we have
$$
 K^{[k]} m=  a^{[k]}(r_k^{-1})=m r_k^{-m\tau+\tau}.
$$
All these identities are satisfied for $\tau \in \RR$ and $K^{[k]}=r_k^{-m\tau+\tau}$.
\eop

In the following remark we list several important properties of the exponential
B-splines.

\medskip \noindent
\begin{remark} \label{rem:exp_b_splines_generation}
$(i)$ The scheme associated with
\begin{equation}\label{symbol:exp-pol_particular}
  a^{[k]}(z)=
   \left(\sum_{\varepsilon \in E} r_k^\varepsilon z^\varepsilon\right)
    \left( \sum_{\varepsilon \in E} s_k^\varepsilon z^\varepsilon \right)
\end{equation}
for $r_k=e^{\lambda m^{-(k+1)}}$ and $s_k=e^{\mu m^{-(k+1)}}$, $\lambda, \mu \in \CC$, $k \ge 0$,
is convergent, \cite{DL95, DL}. It has the basic limit function $\phi=\phi_1 * \phi_2$ with
$\phi_1(x):=e^{\lambda  x}\chi_{[0,1)}$ and $\phi_2(x)=e^{\mu x}\chi_{[0,1)}$.
The function $\phi_1 * \phi_2$ is $C^0$, locally
an exponential function on $[0,1)$ and $[1,2)$, globally supported on $[0,2)$ and is a linear
combination of $e^{\lambda  x}$ and $e^{\mu x}$. 

\noindent $(ii)$
In general, the $n$-fold convolution
$
\phi=\beta_1 \ast \beta_2\ast \cdots \ast \beta_n
$
is  $C^{n-2}$, locally an exponential function on $[J-1,J),\ J=1,\ldots,n$, globally supported on
$[0,n)$ and is a linear combination of the corresponding exponential functions. Such a function
$\phi$  is a basic limit function of the non-stationary scheme given by $n$-fold products of
the symbols in \eqref{symbol:exp-poluni}. Thus, exponential polynomials   are generated by
non-stationary subdivision schemes
with factors in \eqref{symbol:exp-poluni}. 

\noindent $(iii)$ By \cite[Theorem 4.3]{Ron92} and the compact support of the masks,
the exponential B-splines satisfy the assumptions of Proposition
\ref{prop:independence_vs_nonsingular}, if no two purely imaginary $\lambda$ and
$\mu$   satisfy $m^{-\ell}(\lambda-\mu)=2\pi i \ell$ for $\ell \ge L$, $L \ge 0$. Thus, the corresponding non-stationary schemes are
non-singular.
\end{remark}

\smallskip \noindent In the case of several exponential factors the space $EP_{\Gamma, \Lambda}$
that is reproduced by the associated non-stationary subdivision scheme is at most
of dimension $2$. The results of Propositions \ref{prop:Yoon1} and \ref{prop:Yoon2} are consistent with the observation in
\cite{Yoon} for the case $m=2$.

\begin{proposition} \label{prop:Yoon1}
Let $\Gamma=\{0\}$ and $\Lambda=\{\lambda\}$, $\lambda \in \CC$.
Then the non-singular non-stationary scheme defined by
\begin{equation}\label{eq:simbol1f1}
    a^{[k]}(z)=K^{[k]} (1+r_k z+r^2_k z^2+\cdots+r^{m-1}_k z^{m-1})^n, \  r_k=e^{\lambda
 m^{-(k+1)}},\  n \ge 2,
\end{equation}
reproduces at most $EP_{\Gamma, \Lambda}=\{e^{\lambda x}, xe^{\lambda x}\}$,
if $K^{[k]}=m^{1-n}r_k^{-m\tau+\tau}$ and $\tau \in \RR$.
\end{proposition}
\pf
We only need to check the  algebraic conditions  $a^{[k]}(v)=mv^{m\tau-\tau}$ and  $ v \frac{da^{[k]}(z)}{dz}_{z=v}=mv^{m\tau -\tau}\tau$
at $v=r_k^{-1}=e^{-\lambda\,m^{-(k+1)}}$ as the rest of the conditions in \eqref{algebraic_conditions}
are trivially satisfied.
By Theorem \ref{th:characterization_poly_reproduction} the scheme reproduces $EP_{\Gamma, \Lambda}$
if and only if
\begin{eqnarray*}
  K^{[k]} m^n &=& m r_k^{-m\tau+\tau}\, \\
 r_k^{-1} K^{[k]} n r_k m^{n-1} \frac{m(m-1)}{2} &=& m r_k^{-m\tau+\tau} \tau\,.
 \end{eqnarray*}
The first two identities imply that $K^{[k]} = m^{1-n} r_k^{-m\tau+\tau}$ with
$\tau=\frac{n(m-1)}{2}$. Next, to guarantee the reproduction of the exponential polynomial $x^2e^{\lambda x}$,
the algebraic condition $ v^2 \frac{d^2a^{[k]}(z)}{dz^2}_{z=v}=mv^{m\tau-\tau}\tau(\tau-1)$ should be
satisfied. Or, equivalently, we get
$$
 r_k^{-2} K^{[k]} \left( n(n-1) r_k^{2} m^{n-2} \frac{(m-1)^2m^2}{4}+
  n\, m^{n-1}r_k^{2}\frac13m(m-1)(m-2)\right) =m r_k^{-m\tau+\tau} \tau(\tau-1),
$$
that becomes
$$
 \frac{1}{4}n(n-1)(m-1)^2+\frac{1}{3}n(m-1)(m-2)= \frac{n(m-1)}{2} \left( \frac{n(m-1)}{2}-1\right),
$$
or, equivalently,
$$
 3(n-1)(m-1)+4(m-2)=3n(m-1)-6
$$
 which can be only satisfied for $m=-1$.
\eop

\smallskip \noindent We continue with the analysis of another case.

\begin{proposition} \label{prop:Yoon2}
Let $\Gamma=\{0\}$ and $\Lambda=\{\lambda, \mu\}$, $\lambda, \mu \in \CC$,
$\lambda \not =\mu$.
Then the non-singular non-stationary scheme given by
\begin{equation} \label{eq:two_n}
    a^{[k]}(z)= K^{[k]}
   \left(\sum_{\varepsilon \in E} r_k^\varepsilon z^\varepsilon \right)^n
   \left( \sum_{\varepsilon \in E} s_k^\varepsilon z^\varepsilon \right)^n, \quad r_k=e^{\lambda m^{-(k+1)}},
    \quad  s_k=e^{\mu m^{-(k+1)}}
\end{equation}
reproduces at most $EP_{\Gamma, \Lambda}=\{e^{\lambda x}, e^{\mu x}\}$,
if
$$
 K^{[k]}=m^{1-n}  \left(\sum_{\varepsilon \in E}  r_k^{m-1-\varepsilon} s_k^\varepsilon\right)^{-n}
 \quad \hbox{and}\quad \tau=n.
$$
\end{proposition}
\pf
We only need to ensure that $a^{[k]}(v)=mv^{m\tau-\tau}$ for $v \in \{r^{-1}_k,\ s^{-1}_k\}$
as the rest of the conditions in \eqref{algebraic_conditions} are trivially satisfied.
By Theorem \ref{th:characterization_poly_reproduction}, the scheme reproduces $EP_{\Gamma, \Lambda}$
if and only if
\begin{eqnarray*}
 K^{[k]}\, m^n \left(\sum_{\varepsilon \in E}  r_k^{-\varepsilon} s_k^\varepsilon\right)^n&=&
        m\, r_k^{-m\tau+\tau},\\
  K^{[k]}\, m^n \left(\sum_{\varepsilon \in E}  r_k^{\varepsilon} s_k^{-\varepsilon}\right)^n &=&
        m s_k^{-m\tau+\tau}.
\end{eqnarray*}
Therefore, the reproduction of $EP_{\Gamma, \Lambda}$ is possible if and only if
\begin{eqnarray*}
 K^{[k]}\, m^n r_k^{n(1-m)} \left(\sum_{\varepsilon \in E}  r_k^{m-1-\varepsilon} s_k^\varepsilon\right)^n&=&
        m\, r_k^{-m\tau+\tau} \\
  K^{[k]}\, m^n s_k^{n(1-m)} \left(\sum_{\varepsilon \in E}  r_k^{\varepsilon} s_k^{m-1-\varepsilon}\right)^n &=&
        m s_k^{-m\tau+\tau}\,,
\end{eqnarray*}
or, equivalently,
$$
 r_k^{(\tau-n)(1-m)}=s_k^{(\tau-n)(1-m)},
$$
which is the case if and only if either $\tau=n$ or $\lambda=\mu$. For $\tau=n$ the conditions on the
first derivative of $a^{[k]}(z)$ at $z=r_k^{-1}$ are satisfied if and only if $\lambda=\mu$.
\eop

\begin{remark}  \label{rem:Yoon3}
Increasing the number of different factors in \eqref{eq:two_n} to $3$ would lead to the system of equations
\begin{eqnarray*}
   \left(\sum_{\varepsilon \in E}  w_k^{m-1-\varepsilon} s_k^\varepsilon\right)^n  r_{k}^{(\tau-2n)(1-m)}
  &=& \left(\sum_{\varepsilon \in E}  w_k^{m-1-\varepsilon} r_k^\varepsilon\right)^n  s_{k}^{(\tau-2n)(1-m)} \\&=&
   \left(\sum_{\varepsilon \in E}  r_k^{m-1-\varepsilon} s_k^\varepsilon\right)^n w_{k}^{(\tau-2n)(1-m)}
\end{eqnarray*}
in four unknowns $\tau$, $r_k=e^{\lambda m^{-(k+1)}}$, $s_k=e^{\mu m^{-(k+1)}}$ and $w_k=e^{\eta m^{-(k+1)}}$.
Numerical experiments show that these are satisfied when either $\lambda=\eta=\mu$ or
only two of $\lambda,\eta,\mu$ are different. The analysis of other possible
schemes is not one of our goals.
\end{remark}

\subsection{Exponential box splines} \label{subsec:expboxsplines}

In this subsection we describe the structure of the
symbols of non-stationary schemes associated with exponential box splines. We
use them in subsection \ref{subsec:butterfly} to define a non-stationary butterfly scheme.

\smallskip
Let $M=nI$, $n \ge 2$. Then $m=|\hbox{det}M|=n^s$ and $E=\{0,n-1\}^s$.

\begin{proposition} \label{prop:box_spline}
Let
\begin{equation}\label{symbol:exp-pol}
  a^{[k]}(\bz)= \sum_{\varepsilon \in E } \br_k^{\varepsilon} \bz^\varepsilon,   \quad
  k \ge 0, \quad \br_k=e^{\blambda \cdot M^{-k-1}}, \quad  \blambda \in \CC^s \,.
\end{equation}
Then the basic limit function
is $\phi(\bx)=e^{\blambda \cdot \bx} \chi_{[0,1)^s}$.
\end{proposition}
\pf
As the symbols $ a^{[k]}(\bz)$ have a tensor-product structure,
the proof is a straightforward  generalization of our univariate result on exponential B-splines. We only need to observe that (see e.g. \cite{GH}), even for general dilation matrices, every point $\bx$ in the support of
$\phi$ has the representation
\begin{equation}\label{def:x_M_adic}
 \bx= \sum_{j=1}^\infty M^{-j} \bvarepsilon_j, \quad \bvarepsilon_j \in E\,.
\end{equation}
The difficulty of studying the continuity of $\phi$ arises only at the points having at least one $M-$adic component, i.e.,
$$
   x_\ell= \sum_{j=1}^k \left( M^{-j} \bvarepsilon_j \right)_\ell =
   \sum_{j=1}^{k-1} \left(  M^{-j} \bvarepsilon_j \right)_\ell +
    n^{-k} (\bvarepsilon_{k,\ell}-1) +  \sum_{j=k+1}^{\infty} M^{-j} (n-1) \bone_\ell, \  \bvarepsilon_{k,\ell} \not=0\,,
$$
where $\bone_\ell$ is the standard $\ell-$th unit vector of $\RR^s$ and
\begin{eqnarray*}
 \sum_{j=k+1}^\infty M^{-j} (n-1) \bone&=&M^{-k-1} \left( \sum_{j=0}^\infty M^{-j} \right)  (n-1) \bone \\
 &=&  M^{-k-1} \left( I-M^{-1} \right)^{-1}  (n-1) \bone\\ &=& M^{-k}  \left( M-I \right)^{-1}  (n-1)
 \bone=M^{-k} \bone\,.
\end{eqnarray*}
\eop

\begin{remark} \label{rem:box_splines}
By \cite[(4.1)]{Ron_exp_box_splines} or by \cite[example 6]{DL95}, the scheme defined by the
masks in \eqref{symbol:exp-pol} converges weakly to an exponential
box spline $e^{\blambda \bx} \chi_{[0,1)^s}$, where $\chi_{[0,1)^s}$ is
a characteristic function of $[0,1)^s$. To get the similar result for $\phi$ with sets of directions
$\Upsilon=\{\upsilon_1, \dots \upsilon_s\}$,
$\hbox{span}\Upsilon=\RR^s$, other than $E$, one just carries out
the affine transformation $\by=[\upsilon_1 \dots \upsilon_s] \bx$
for $\bx \in [0,1)^s$. See \cite[example 6]{DL95} for the masks of the smoother exponential
box splines.
\end{remark}

\subsection{A $C^2$ binary, non-stationary, dual 4-point subdivision
scheme reproducing conics} \label{subsec:examples_uni1}

In this and the following subsections we show that our algebraic conditions in \eqref{algebraic_conditions} can be
efficiently used for constructing new univariate subdivision schemes with desired reproduction
properties and with enhanced smoothness.

We start by deriving the non-stationary counterpart of the binary dual
4-point subdivision scheme in \cite{DFH05}. We show that, compared to
the binary non-stationary interpolatory 4-point scheme in
\cite{BCR07a}, our scheme reproduces the same space of
exponential polynomials, \ie \, $span \{1,x,e^{\lambda x},e^{-\lambda x}\}$ with 
$\lambda \in \left(\RR \cup i \RR\right) \setminus \{0\}$, but it is $C^2$ instead of $C^1$.
Furthermore, among all existing non-stationary binary schemes that are $C^2$ and reproduce
conic sections, see \cite{CGR2011,CR2010,CR2011,DLL03,Yoon,JLYoon,R09}, our scheme turns out to be the one with the smallest support.

\smallskip \noindent
In order to define our scheme, we consider a sequence
of masks of the form
$$
\ba^{[k]}=\{\cdots,\, 0, \, c^{[k]}_{3,2}, \, c^{[k]}_{0,1}, \, c^{[k]}_{2,2}, \,
c^{[k]}_{1,1}, \, c^{[k]}_{1,2}, \, c^{[k]}_{2,1}, \, c^{[k]}_{0,2}, \, c^{[k]}_{3,1},\,
0\, \cdots\},\quad k\ge 0,
$$
and  make use of Corollary  \ref{cor:multi_repr} to determine $c^{[k]}_{i,j}$,
$i=0,\ldots,3$, $j=1,2$, and the shift parameter $\tau \in \RR$ such that the space $span \{1,x,e^{\lambda x},e^{-\lambda x}\}$ is reproduced.
Namely, let $\Lambda=\{0,\lambda,-\lambda\}$, $\lambda \in \left(\RR \cup i \RR\right) \setminus \{0\} $, and $\Gamma=\{0,1\}$.
Then, for $E=\{0,1\}$ and $k \ge 0$,  we have
$$
V_k=\left \{ e^{\pi\, i \varepsilon} e^{-\lambda \,2^{-(k+1)}}\ : \ \varepsilon \in E, \ \lambda \in \Lambda \right \}
$$
and
$$
V'_k=\left \{ -e^{-\lambda \,2^{-(k+1)}}\ : \ \lambda \in \Lambda \right \}=\{-e^{\lambda 2^{-(k+1)}}, -e^{-\lambda 2^{-(k+1)}}, -1 \}.
$$
Thus, by Corollary  \ref{cor:multi_repr} and Remark \ref{rem:Q} $(ii)$,
 for $Q=\{(0,0),(1,0), (0,\lambda), (0, -\lambda)\}$, we need to solve the following linear system of $8$ equations in $8$ unknowns
$$
\left \{
\begin{array}{l}
a^{[k]}(v)=0, \quad \forall v \in V_k',\smallskip\\
a^{[k]}(v)=2\,v^{\tau},\quad \forall v \in \{e^{\lambda 2^{-(k+1)}}, e^{-\lambda 2^{-(k+1)}}, 1 \}, \smallskip\\
\begin{array}{ll}
\frac{da^{[k]}(z)}{dz}|_{z=v}=0, & \hbox{for} \ v=-1, \smallskip\\
\frac{da^{[k]}(z)}{dz}|_{z=v}=2 \tau, & \hbox{for} \ v=1.
\end{array}
\end{array}
\right.
$$
We get $\tau=-\frac12$ and
\begin{eqnarray*}
c^{[k]}_{0,2}=c^{[k]}_{0,1}&=&-\frac{6(w^{[k]})^2 + 2w^{[k]} - 1}{64(w^{[k]})^3
(2(w^{[k]})^2 - 1)(w^{[k]} + 1)}, \nonumber \smallskip\\
c^{[k]}_{1,2}=c^{[k]}_{1,1}&=&\frac{10(w^{[k]})^2 + 2w^{[k]} - 3}{64(w^{[k]})^3
(2(w^{[k]})^2 - 1)(w^{[k]} + 1)} + \frac{3}{4}, \smallskip\\
c^{[k]}_{2,2}=c^{[k]}_{2,1}&=&\frac{-2(w^{[k]})^2 + 2w^{[k]} +
3}{64(w^{[k]})^3(2(w^{[k]})^2 - 1)(w^{[k]} + 1)} + \frac{1}{4},
\nonumber \smallskip\\
c^{[k]}_{3,2}=c^{[k]}_{3,1}&=&-\frac{2(w^{[k]})^2 + 2w^{[k]} +
1}{64(w^{[k]})^3(2(w^{[k]})^2 - 1)(w^{[k]} + 1)}, \nonumber
\end{eqnarray*}
with
$$
w^{[k]}=\frac{1}{2}(e^{2^{-(k+1)}\lambda /2}+e^{-2^{-(k+1)}\lambda/2})\,.
$$
In conclusion, the non-stationary dual 4-point subdivision scheme we propose is defined
by the $k$-level symbol
$$
\begin{array}{l}
a^{[k]}(z)=-z^{-4} \, \frac{1}{64(w^{[k]})^3(2(w^{[k]})^2-1)(w^{[k]}+1)} \, (z +
1)^3 \, \Big(z^2+(4(w^{[k]})^2-2)z+1\Big) \, \cdot\\
\cdot  \, \Big((2(w^{[k]})^2 + 2w^{[k]} + 1)z^2
-(8(w^{[k]})^4+8(w^{[k]})^3+2)z + 2(w^{[k]})^2 + 2w^{[k]} + 1\Big).
\end{array}
$$
We observe that
$$
\lim_{k\rightarrow \infty} c^{[k]}_{0,1}=-\frac{7}{128}, \quad
\lim_{k\rightarrow \infty} c^{[k]}_{1,1}=\frac{105}{128}, \quad
\lim_{k\rightarrow \infty} c^{[k]}_{2,1}=\frac{35}{128}, \quad
\lim_{k\rightarrow \infty} c^{[k]}_{3,1}=-\frac{5}{128},
$$
i.e. when $k$ tends to infinity $\ba^{[k]}$ converges to the mask of the above mentioned
stationary dual 4-point subdivision scheme. More precisely, the result in \cite[Theorem 8]{DL95},
implies that our non-stationary scheme is \emph{asymptotically equivalent} (of ``order'' 2) to the stationary dual
4-point subdivision scheme in \cite{DFH05}. This property allows us to conclude that our scheme is indeed $C^2$.

\subsection{A $C^2$ ternary, non-stationary, dual 4-point
subdivision scheme reproducing conics} \label{subsec:examples_uni2}
In this subsection we derive  a $C^2$ ternary, non-stationary,
dual 4-point subdivision scheme reproducing the space $span \{1,x,e^{\lambda x},e^{-\lambda x}\}$ with $\lambda \in \left(\RR \cup i \RR\right) \setminus \{0\}$.
 Compared with the ternary, non-stationary, interpolatory 4-point scheme in \cite{BCR07b}, our scheme will have an additional feature of reproducing conic sections. Compared with the ternary, non-stationary, interpolatory 4-point scheme that reproduces the same space of exponential polynomials
 and is $C^1$ \cite{BCR09,BCR10},  our scheme will be $C^2$.\\
As in the previous subsection,  we start by defining a sequence of masks of the form
$$\ba^{[k]}=\{\cdots,\, 0, \, c^{[k]}_{0,3}, \, c^{[k]}_{0,2}, \,
c^{[k]}_{0,1}, \, c^{[k]}_{1,3}, \, c^{[k]}_{1,2}, \, c^{[k]}_{1,1},
\, c^{[k]}_{2,3}, \, c^{[k]}_{2,2}, \, c^{[k]}_{2,1}, \,
c^{[k]}_{3,3}, \, c^{[k]}_{3,2}, \, c^{[k]}_{3,1}, \, 0, \, \cdots\}\,,
$$
and we derive the coefficients $c^{[k]}_{i,j},\ i=0,\ldots,3$, $j=1,2,3$, and the shift parameter
$\tau \in \RR$ by requiring that the scheme reproduces $span \{1,x,e^{\lambda x},e^{-\lambda x}\}$. Let $\Lambda=\{0,\lambda,-\lambda\}$,
$\lambda \in \left(\RR \cup i \RR\right) \setminus \{0\}$, and
$\Gamma=\{0,1\}$. Then, for $E=\{0,1,2\}$ and $k \ge 0$, we have
$$
V_k=\left \{ e^{ \frac{2\pi\, i\varepsilon}{3}} e^{-\lambda\,3^{-(k+1)}}\ : \ \varepsilon\in E, \ \lambda \in \Lambda \right\}
$$
and
$$
V_k'=\left \{ e^{ \frac{2\pi\, i\varepsilon}{3}} e^{-\lambda\,3^{-(k+1)}}\ : \ \varepsilon\in E \setminus \{0\}, \
\lambda \in \Lambda \right \}.
$$
The corresponding algebraic conditions of Corollary  \ref{cor:multi_repr} and Remark \ref{rem:Q} $(ii)$
with $Q=\{(0,0),(1,0), (0,\lambda), (0, -\lambda)\}$ lead to the following system of $12$ equations in $12$ unknown
$$
\left \{
\begin{array}{l}
a^{[k]}(v)=0, \quad \forall v \in V_k',\\
a^{[k]}(v)=3\,v^{2\tau},\quad \forall v \in \{e^{\lambda\,3^{-(k+1)}}, e^{-\lambda\,3^{-(k+1)}}, 1 \}, \\
\begin{array}{ll}
\frac{da^{[k]}(z)}{dz}|_{z=v}=0, \quad v=e^{ \frac{2\pi\, i\varepsilon}{3}}, \quad
\varepsilon=1,2, \smallskip\\
\frac{da^{[k]}(z)}{dz}|_{z=v}=6 \tau \quad \hbox{for} \ v=1.
\end{array}
\end{array}
\right.
$$
Solving the above system we get $\tau=-\frac14$,
\begin{eqnarray*}
c^{[k]}_{0,2}&=&-\frac{1}{8w^{[k]}(2w^{[k]} - 1)^2(4(w^{[k]})^2 - 3)(w^{[k]} + 1)}, \smallskip
\nonumber \\
c^{[k]}_{1,2}&=&\frac{1}{8w^{[k]}(2w^{[k]} - 1)^2(4(w^{[k]})^2 - 3)(w^{[k]} + 1)} + \frac12,
\smallskip  \\
c^{[k]}_{2,2}&=&\frac{1}{8w^{[k]}(2w^{[k]} - 1)^2(4(w^{[k]})^2 - 3)(w^{[k]} + 1)} + \frac12,
\smallskip \nonumber \\
c^{[k]}_{3,2}&=&-\frac{1}{8w^{[k]}(2w^{[k]} - 1)^2(4(w^{[k]})^2 - 3)(w^{[k]} + 1)} \nonumber
\end{eqnarray*}
and
\begin{eqnarray*}
c^{[k]}_{0,3}=c^{[k]}_{3,1} &=&\frac{16(w^{[k]})^4 + 16(w^{[k]})^3 + 3}{24w^{[k]}(4(w^{[k]})^2 -
1)^3(- 4(w^{[k]})^3 - 4(w^{[k]})^2 + 3w^{[k]} + 3)}, \smallskip \nonumber\\
c^{[k]}_{1,3}=c^{[k]}_{2,1} &=&-\frac{16(w^{[k]})^4 - 16(w^{[k]})^2 - 4w^{[k]} - 1}{8w^{[k]}(2w^{[k]} -
1)^3(2w^{[k]} + 1)^3(4(w^{[k]})^2 - 3)(w^{[k]} + 1)} + \frac16, \smallskip \nonumber\\
c^{[k]}_{2,3}=c^{[k]}_{1,1} &=&\frac{48(w^{[k]})^4 + 16(w^{[k]})^3 - 32(w^{[k]})^2 - 8w^{[k]} +
1}{8w^{[k]}(2w^{[k]} - 1)^3(2w^{[k]} + 1)^3(4(w^{[k]})^2 - 3)(w^{[k]} + 1)} + \frac56, \smallskip
\nonumber \\
c^{[k]}_{3,3}=c^{[k]}_{0,1} &=&\frac{80(w^{[k]})^4 + 32(w^{[k]})^3 - 48(w^{[k]})^2 - 12w^{[k]} +
3}{24w^{[k]}(4(w^{[k]})^2 - 1)^3(- 4(w^{[k]})^3 - 4(w^{[k]})^2 + 3w^{[k]} + 3)}
\end{eqnarray*}
where
$$
w^{[k]}=\frac{1}{2}\left( e^{3^{-(k+1)}\lambda /2}+e^{-3^{-(k+1)}\lambda /2} \right).
$$
The resulting $k$-level symbol is
$$
\begin{array}{l}
a^{[k]}(z)=-z^{-6} K^{[k]} \  (z^2 + z + 1)^2 (z + 1)  \\
\Big( z^4 + (4(w^{[k]})^2 - 2)z^3 + (16(w^{[k]})^4 - 16(w^{[k]})^2 + 3)z^2 + (4(w^{[k]})^2 - 2)z + 1 \Big) \cdot \\
\Big( (16(w^{[k]})^4 + 16(w^{[k]})^3 + 3)z^2 + (- 64(w^{[k]})^6 - 64(w^{[k]})^5 + 32(w^{[k]})^4 + \\
\quad 32(w^{[k]})^3 - 12(w^{[k]})^2 - 12w^{[k]} - 6)z + 16(w^{[k]})^4 + 16(w^{[k]})^3 + 3 \Big),
\end{array}
$$
where
$$
 K^{[k]}= \frac{1}{24 w^{[k]} (2w^{[k]} - 1)^3 (2w^{[k]} + 1)^3 (4(w^{[k]})^2 - 3) (w^{[k]} + 1)}.
$$
This non-stationary dual 4-point scheme is ``order'' 2 asymptotically equivalent to the stationary
dual 4-point scheme with the mask
$$
\lim_{k\rightarrow \infty} \ba^{[k]}=\left \{ -\frac{35}{1296}, \,
-\frac{1}{16}, \, -\frac{55}{1296}, \, \frac{77}{432}, \,
\frac{9}{16}, \,
\frac{385}{432}, \, \frac{385}{432}, \, \frac{9}{16}, \,
\frac{77}{432}, \, -\frac{55}{1296}, \, -\frac{1}{16}, \,
-\frac{35}{1296}  \right \}.
$$
and its symbol satisfies
$$
\lim_{k\rightarrow \infty}  a^{[k]}(z)=-z^{-6}\frac{1}{1296}(z^2 + z + 1)^4(z + 1)(35z^2 - 94z + 35).
$$
By \cite[Section 5.2]{ContiRomani13}, this stationary scheme is $C^2$. Thus, applying \cite[Theorem 8]{DL95}
we deduce that our non-stationary scheme is also $C^2$.

\medskip
The rest of this section is devoted to the multivariate case.

\subsection{Non-stationary butterfly scheme} \label{subsec:butterfly}
In this subsection, we derive a non-stationary butterfly scheme using the results in \cite{CCJZ}.
Let $\br_k=e^{\blambda 2^{-(k+1)}}$, $\blambda \in \RR^2 \cup i\RR^2$ and $k \ge 0$. The symbols of the
non-stationary butterfly scheme we construct are combinations of the symbols
$$
 B^{[k]}_{j,h,\ell}(z_1,z_2)=\left(\frac{1+ r_{k,1} z_1}{2}
 \right)^j \left(\frac{1+ r_{k,2} z_2}{2}
 \right)^h \left(\frac{1+ r_{k,1} r_{k,2} z_1z_2}{2}
 \right)^\ell, \ j,h,\ell \in \NN,
$$
corresponding to the three-directional exponential box splines. In this case, $M=2I$,
$$
 E=\left\{ (0,0)^T, (1,0)^T, (0,1)^T, (1,1)^T \right\}
$$
and
$$
 \Xi=\left\{ (1,1)^T, (-1,1)^T, (1,-1)^T, (-1,-1)^T \right\}.
$$
It is easy to check that the symbols
\begin{eqnarray*}
 a^{[k]}(\bz)=4\, \Big( 7 \, r_{k,1} r_{k,2} z_1z_2 B^{[k]}_{2,2,2}(\bz)
    &-& 2\, r_{k,1} z_1 B^{[k]}_{1,3,3}(\bz)
     -  2\, r_{k,2} z_2 B^{[k]}_{3,1,3}(\bz) \\
     &-& 2\, r_{k,1} r_{k,2} z_1z_2 B^{[k]}_{3,3,1}(\bz) \Big),
    \quad k\ge 0,
\end{eqnarray*}
satisfy \eqref{algebraic_conditions_multi_cor} for
$\Gamma=\{\bgamma \in \NN_0^2\ : \  |\bgamma| < 4\}$ and for
$$
 V'_k= \left\{ (-r^{-1}_{k,1}, r^{-1}_{k,2})^T, (r^{-1}_{k,1},-r^{-1}_{k,2})^T, (-r^{-1}_{k,1},-r^{-1}_{k,2})^T
 \right\}.
$$
Thus, the scheme  $S_{\{\ba^{[k]}, \ k\ge 0\}}$ is
$EP_{\blambda, \Gamma}-$generating.  The scheme is also
interpolatory, since the symbols satisfy
$$
 z_1^{-3} z_2^{-3} \left(
 a^{[k]}(z_1,z_2)+a^{[k]}(-z_1,z_2)+a^{[k]}(z_1,-z_2)+a^{[k]}(-z_1,-z_2)\right)=4, \ k \ge 0.
$$
Thus, by Corollary \ref{cor:interp}, it is also $EP_{\blambda, \Gamma}-$reproducing with
$\btau=(0,0)^T$. We observe that our algebraic conditions allow us to simplify the
construction of the non-stationary butterfly scheme in \cite{Yoon1} and to improve its reproduction properties.
The smoothness analysis
of this scheme is done using standard techniques in \cite{DL}.

\subsection{Non-stationary schemes obtained by convolution}
Another example is of a non-stationary subdivision scheme
with dilation matrix $ M=\left( \begin{array}{cc} 2&1\\0&2\end{array}\right)$ with $m=4$.
In this case we have
$$
 E=\left\{ (0,0)^T, (1,0)^T, (1,1)^T, (2,1)^T \right\}, \quad
 \Xi=\left\{ (1,1)^T, (1,-1)^T,  (-1,i)^T, (-1,-i)^T   \right\},
$$
since the coset representatives of $\ZZ^s / M^T \ZZ^s$ are $\{(0,0)^T, (1,0)^T, (1,1)^T, (1,2)^T\}$.
The results in \cite{CDM,Dahlke} imply that the stationary scheme associated with the symbol
$$
  a(\bz)=\frac{1}{4}b^2(\bz), \quad  b(\bz)=\sum_{\bvarepsilon \in E} \bz^\bvarepsilon,
$$
is convergent. Thus, by \cite{DL95}, the non-stationary scheme given by the symbols
$$
 a^{[k]}(\bz)=\frac{1}{4}\left(b^{[k]}(\bz)\right)^2, \quad  b^{[k]}(\bz)=\sum_{\bvarepsilon \in E}
  \br_k^\bvarepsilon \bz^\bvarepsilon, \quad k \ge 0,
$$
is also convergent, if $\br_k$ satisfy $\displaystyle \sum_{k \in
\ZZ_+} |\bone-\br_k| < \infty$. For example, we let
$\br_k=e^{ \blambda \cdot M^{-(k+1)}}$ for some $\blambda \in \RR^2 \cup i \RR^2$. The
conditions in \eqref{algebraic_conditions_multiGENE} are satisfied for the associated $V_k'$. Thus,
this non-stationary scheme generates polynomials  $\{ \bx^\bgamma e^{\blambda \cdot \bx} \ : \
 \bgamma \in \NN_0^s, \ |\bgamma|\le 1 \}$. Furthermore,
due to $q_{(0,0)}(\bz)=1$, we get that the condition
$$
 4=a^{[k]}(\bv)=4 \cdot \bv^{M\btau-\btau}, \quad
 \bv=(r^{-1}_{k,1},r^{-1}_{k,2}),
$$
is satisfied if and only if $\btau=(0,0)^T$. Moreover, e.g. for
$\bgamma=(1,0)^T$, the condition
$$
 8= \bv^{(1,0)} D^{(1,0)} a^{[k]}(\bv)=4 \cdot \bv^{M\btau-\btau}
 q_{(1,0)}(\bnull)= 4 \cdot \bv^{M\btau-\btau}
$$
cannot be satisfied for $\bv=(r^{-1}_{k,1},r^{-1}_{k,2})$ and
$\btau=(0,0)^T$. Thus, we can only conclude that the non-stationary scheme reproduces
$e^{\blambda \cdot \bx}$ with this shift parameter $\btau$.

We show next how to rescale the symbols $a^{[k]}(\bz)$ appropriately in order to enlarge the space 
of exponential polynomials that is being reproduced. To do so, let
$a^{[k]}(\bz)=K^{[k]}\left(b^{[k]}(\bz)\right)^2$, where $K^{[k]}$, $k \ge 0$, are to be
determined from our algebraic conditions. For $\bv=(r_{k,1}^{-1},
r_{k,2}^{-1})$, we get the system of equations
\begin{eqnarray*}
 && 4 K^{[k]}=\bv^{M\btau-\btau}, \quad M\btau-\btau=(\tau_1+\tau_2,\tau_2)^T,\\
 && 8 K^{[k]}=\bv^{M\btau-\btau} q_{(1,0)}(M\btau-\btau)=\bv^{M\btau-\btau}(\tau_1+\tau_2),\\
 && 4 K^{[k]}=\bv^{M\btau-\btau} q_{(0,1)}(M\btau-\btau)=\bv^{M\btau-\btau}\tau_2.
\end{eqnarray*}
The unique solution of this system is
$$
  K^{[k]}=\frac{1}{4} \bv^{M \btau-\btau} \quad  \hbox{and} \quad \btau=(1,1)^T.
$$
Thus, the non-stationary scheme given by the symbols
$$
 a^{[k]}(\bz)=\frac{1}{4} \br_k^{-M\btau+\btau} \left(b^{[k]}(\bz)\right)^2, \quad  b^{[k]}(\bz)=\sum_{\bvarepsilon \in E}
  \br_k^\bvarepsilon \bz^\bvarepsilon, \quad k \ge 0,
$$
generates and reproduces the space $\{\bx^\bgamma e^{\blambda
\cdot \bx}\ :  \ \bgamma \in \NN_0^s, \ |\bgamma| \le 1\}$ for
some $\blambda \in \RR^2 \cup i \RR^2$ and $\btau=(1,1)^T$.

\subsection{$\sqrt{3}-$ subdivision}

The results of this paper also generalize the results of \cite{CC2012} on polynomial reproduction of stationary schemes to the case
of a general dilation matrix. 
This example shows how to determine the degree of polynomial reproduction of the approximating
$\sqrt{3}-$subdivision schemes given in \cite[page 21]{JiangOswald}
from the corresponding mask symbol $a(\bz)$ instead of the iterated symbol $a(z_1z_2^{-2},z_1^2z_2^{-1}) \cdot a(\bz)$, as it is done
in \cite{CC2012}. We also show how to use affine combinations of these schemes to improve their
degree of polynomial reproduction, i.e. corresponding to $\Lambda=\{\bnull\}$. 

The dilation matrix in this case is
$$
  M=\left( \begin{array}{rr} 1&2\\-2&-1\end{array}\right), \quad M^2=-3I,
$$
and the mask symbol is given by
\begin{eqnarray}
 a(\bz)&=&\frac{1}{6}\left(z_1z_2+z_1^{-1}z_2^{-1}+z_1^{-1}z_2^2+z_1^{-2}z_2+z_1z_2^{-2}+z_1^2z_2^{-1}\right)+ \notag \\
  &&\frac{1}{3} \left(z_1^{-1}+z_2+z_1z_2^{-1}\right)+\frac{1}{3}
 \left(z_2^{-1}+z_1+z_1^{-1}z_2\right). \notag
\end{eqnarray}
The associated subdivision scheme satisfies zero conditions of at most
order $2$, see \cite{JiangOswald}. The result of Corollary \ref{cor:multi_repr} 
yields $\btau=(0,0)$ in \eqref{def:parametrization}, which implies
that the corresponding scheme reproduces linear polynomials. Since, the mask
symbol satisfies at most zero conditions of order $2$, the associated
refinable function has approximation order $2$, see \cite{ALevin}. Note that,
similarly, the corresponding $\btau$ is $(0,0)$ for all
approximating $\sqrt{3}-$subdivision schemes given in \cite[page
21]{JiangOswald}.

Take next the symbols $a_j(\bz)$, $j=1,2,3,4$ of the four approximating subdivision schemes in \cite[page
21]{JiangOswald} that satisfy zero conditions of order $3$ and consider their affine combination
$$
 a(\bz)=\sum_{j=1}^4 \lambda_j \cdot a_j(\bz), \quad \sum_{j=1}^4 \lambda_j=1.
$$
The resulting scheme still satisfies Condition~Z$_3$ and in addition the rest of the conditions in
Corollary \ref{cor:multi_repr}  for $k=3$, if
$$
 \lambda_1=\frac{1}{2} \lambda_3+\lambda_4+2, \quad \lambda_2=-\frac{3}{2} \lambda_3-2\lambda_4-1,\quad
 \lambda_3,\lambda_4 \in \RR.
$$
Surprisingly, for any such choice of the parameters $\lambda_1, \dots, \lambda_4$ the resulting
scheme is given by
$$
 a(\bz)=1-\frac{1}{9}(z_1^{-2}+z_1^{-2}z_2^2+z_2^2+z_1^2+z_1^2z_2^{-2}+z_2^{-2})+
 \frac{4}{9}(z_1^{-1}+z_1^{-1}z_2+z_2+z_1+z_1z_2^{-1}+z_2^{-1})
$$
and defines the interpolatory scheme considered in \cite[page 17]{JiangOswald}. Thus, by \cite{ALevin},
this scheme has approximation order  $3$.

\bigskip
\bibliographystyle{amsplain}


\end{document}